\documentclass[11pt]{article}
\usepackage{hyperref}
\usepackage{breakurl}
\hypersetup{ colorlinks, citecolor=blue, filecolor=black, linkcolor=blue, urlcolor=blue} 
\usepackage[top=1in, bottom=1in, left=1in, right=1in]{geometry}
\usepackage{amsmath}
\usepackage{amsthm}
\usepackage{amssymb}
\usepackage[utf8]{inputenc}
\usepackage[english]{babel}
\usepackage{changepage}
\usepackage{subcaption}

\usepackage{graphicx}

\usepackage{float}
\usepackage{color,soul}
\usepackage{bm}
\numberwithin{equation}{section}
\usepackage[short, c2]{optidef}

\usepackage{mathtools}
\usepackage{interval}
\DeclarePairedDelimiter\abs{\lvert}{\rvert}%
\DeclarePairedDelimiter\norm{\lVert}{\rVert}%
\makeatletter
\let\oldabs\abs
\def\abs{\@ifstar{\oldabs}{\oldabs*}}
\let\oldnorm\norm
\def\norm{\@ifstar{\oldnorm}{\oldnorm*}}
\makeatother

\usepackage{setspace}
\usepackage{fancyhdr}
\usepackage{lipsum}
\usepackage{atbegshi}
\usepackage{enumerate}

\newtheorem{lemma}{Lemma}
\newtheorem{proposition}{Proposition}
\newtheorem{theorem}{Theorem}

\usepackage{algorithm}
\usepackage{algpseudocode,algorithmicx}
\usepackage{booktabs,ragged2e}
\algnewcommand\algorithmicforeach{\textbf{for each}}
\algdef{S}[FOR]{ForEach}[1]{\algorithmicforeach\ #1\ \algorithmicdo}

\usepackage{caption} 
\usepackage{changepage}
\usepackage{array}
\newcolumntype{P}[1]{>{\centering\arraybackslash}p{#1}}
\newcolumntype{M}[1]{>{\RaggedLeft\arraybackslash}m{#1}}
\newcolumntype{C}[1]{>{\centering\arraybackslash}m{#1}}

\usepackage{multirow}

\usepackage{siunitx}
\title{\vspace{-5ex} Consensus-Based Dantzig-Wolfe Decomposition}
\author{
    Mohamed El Tonbari\thanks{Email: {\tt  mtonbari@gatech.edu}} \
    Shabbir Ahmed \\
 School of Industrial \& Systems Engineering,\\ Georgia Institute of Technology   
    }
\date{}

\begin{document}

\maketitle
\raggedbottom 

\begin{abstract}
Dantzig-Wolfe decomposition (DWD) is a classical algorithm for solving large-scale linear programs whose constraint matrix involves a set of independent blocks coupled with a set of linking rows. The algorithm decomposes such a model into a master problem and a set of independent subproblems that can be solved in a distributed manner. In a typical implementation, the master problem is solved centrally. In certain settings, solving the master problem centrally is undesirable or infeasible, such as in the case of decentralized storage of data, or when independent agents who are responsible for the subproblems desire privacy of information. In this paper, we propose a fully distributed DWD algorithm which relies on solving the  master problem using a consensus-based Alternating Direction  Method of Multipliers (ADMM) method. We derive error bounds on the optimality gap and feasibility violation of the proposed approach. We provide preliminary computational results for our algorithm using a Message Passing Interface (MPI) implementation on cutting stock instances from the literature and synthetic instances where we obtain high quality solutions.
\end{abstract}

\section{Introduction}
\subsection{Dantzig-Wolfe Decomposition}
Dantzig-Wolfe decomposition (DWD) \cite{dw1960} is a classical algorithm for solving large-scale linear programs whose constraint matrix involves a set of independent blocks coupled with a set of linking rows. This class of problems are of the form

\begin{mini}
{x}
{\sum_{n=1}^N c_n^\top x_n}
{\label{OPT:P} \tag{$P$}}
{}
\addConstraint{\sum_{n=1}^N A_n x_n}{= t}
\addConstraint{x_n}{\in X_n,}{\quad n = 1, \dots , N}
\end{mini}
where $N$ is the number of blocks, the set $X_n$ denotes the feasible region of the $n$-th block (or subproblem) with the decision vector $x_n$, and $\sum_{n=1}^NA_nx_n = t $ denotes the system of linking constraints. The DWD method decomposes \eqref{OPT:P} into a master problem and a set of $N$ independent subproblems. Throughout this paper, we assume that the sets  $X_n$ are polytopes (i.e. bounded polyhedra) for all $n$. In this case, the master problem is a reformulation where each vector $x_n$ is replaced by a convex combination of the extreme points of $X_n$. From Minkowski's theorem, the master problem can be written as
\begin{mini}
{\lambda}
{\sum_{n=1}^N\sum_{i \in I_n}c_n^\top x_n^i\lambda_{ni}}
{\label{OPT:MP_verbose}}
{}
\addConstraint{\sum_{n=1}^N\sum_{i \in I_n}A_nx_n^i\lambda_{ni}}{= t}
\addConstraint{\sum_{i \in I_n} \lambda_{ni}}{= 1,}{\quad n = 1, \dots, N}
\addConstraint{\lambda_{ni}} {\geq 0,}{\quad \forall i \in I_n, \  n = 1, \dots , N} 
\end{mini}
where $\{x_n^i\}_{i \in I_n}$ and $I_n$ correspond to the extreme points of $X_n$ and their index set, respectively; the variable $\lambda_{ni}$ is the convex multiplier associated with extreme point $i$ of subproblem $n$.

Since the number of extreme points of the subproblem polytopes can be exponentially large, a restricted master problem RMP is solved instead, using only a subset of the extreme points (or columns). Using an optimal dual solution of the RMP,  each block $n$ solves a pricing subproblem independently to find a new extreme point of $X_n$ with negative reduced cost. This process  (called column generation) is repeated until no columns with negative reduced costs can be found. The optimal solution of RMP in the last iteration then provides an optimal solution to \eqref{OPT:P}.

\subsection{Motivation}
DWD has been used to solve a variety of problems, many arising in energy and transportation, and has received great attention in solving mixed-integer programs to obtain tighter relaxation bounds, leading to the popular branch-and-price algorithm. See \cite{bertsimasLP,desaulniers2005,lubbecke2005,nemhauser1988} for details and applications of DWD and column generation.

In each iteration of DWD, the subproblems can be solved independently in a distributed manner but the RMP is solved centrally. 
In certain settings, solving the master problem centrally is undesirable or infeasible. With increasingly large amounts of data available, we are seeing an upward trend in decentralized storage of data, either as a protection against attacks, or due to memory limitations, in which case the data is not available centrally to solve the master problem \cite{boyd}. Alternatively, if we interpret the variables of each block as the decision variables of independent agents, then solving the RMP centrally requires the agents to share data of their constraints, objective, and decisions, potentially violating privacy. Privacy concerns are prevalent in energy, such as in smart-grid optimization, for example. In order to optimize power generation, smart meters provide sensitive and private energy consumption data of the household to the utility provider \cite{zeadally2013towards, belletti2015privacy}. 
Another example is in the healthcare industry with Accountable Care Organizations (ACOs), where a collection of hospitals or clinics form a group to give coordinated care to patients. 
To drive costs down and improve patient care, sensitive patient data would link the optimization problem, causing concerns for privacy violations among ACOs \cite{decamp2014ethical, mcwilliams2016early}.
Developing a consensus-based Dantzig-Wolfe algorithm would permit the use of DWD in such scenarios while handling decentralized storage of data and/or privacy concerns.

\subsection{Contribution}
In this paper, we propose a consensus-based DWD algorithm which relies on solving the dual of the master problem using a consensus-based Alternating Direction  Method of Multipliers (ADMM) algorithm. We address the computational challenges and theoretical questions that arise from solving the master problem in a distributed fashion. To the best of our knowledge, distributed or consensus-based techniques to solve the master problem have not been studied within the context of column generation or Dantzig-Wolfe decomposition.

We illustrate the benefits of working on the dual and using ADMM as a consensus-based method. To reduce computation time, we dynamically adjust the tolerances, where we first aim for inaccurate dual solutions to speed up ADMM convergence, similar to \cite{gondzio1, gondzio3}, and lower the tolerances as the algorithm progresses. This significantly decreases the number of times we need to solve ADMM to high accuracy, yielding computational benefits. 
ADMM is known to be slow to converge to high accuracy, but requires only a few iterations to achieve modest accuracy \cite{boyd}, and thus synergizes well with solving for inaccurate duals at first. Moreover, ADMM has stronger convergence properties than subgradient methods, but performing it on \eqref{OPT:MP_verbose} directly breaks decomposability. We instead solve the dual of \eqref{OPT:MP_verbose}, maintaining decomposability and preserving privacy of information. Indeed, each subproblem need only share their dual variables with a central coordinator, as is common in a privacy preserving setting. By solving the dual and using ADMM, we circumvent primal-recovery issues present in subgradient methods and avoid the need for ergodic sequences \cite{nedic2009approximate, gustavsson2015primal}. We show we can easily recover a primal solution to the original problem that is close to feasible and close to optimal using the Lagrangian multiplers associated with the constraints in the ADMM subproblems. Finally, while most stability techniques in DWD rely on suboptimal duals, our method also handles infeasible duals, as the ADMM approach provides $\epsilon$-optimal and $\delta$-feasible dual solutions.
We prove bounds on the feasibility violation and optimality gap at the recovered primal solution. 

We provide preliminary computational results for the proposed algorithm using a Message Passing Interface (MPI) implementation on cutting stock instances from the literature and synthetic instances where we obtain high quality solutions. Although our main contribution is to tackle decentralized storage of data and privacy, our method also shows a potential computational benefit for instances with a large number of variables and in many of the cutting stock instances.

\subsection{Prior Work}

Solving the dual of RMP in a distributed manner leads to approximate dual solutions used in the pricing subproblems, as opposed to standard DWD where exact optimal dual solutions are readily available. Several stability techniques proposed in the literature use suboptimal but feasible dual solutions to solve the pricing subproblems to circumvent unstable behavior resulting from using optimal dual solutions \cite{gondzio1}. One method is to add a penalty term to reduce the variation of obtained dual solutions \cite{bundle_2006, gondzio2}; another involves using primal-dual interior point method to solve for suboptimal dual solutions which are well-centered, where the optimality tolerance of the interior point method is dynamically adjusted to reduce the computation time of solving the RMP~\cite{gondzio1,gondzio2,gondzio3}. In these methods, suboptimal dual solutions are used to solve the pricing subproblems, where tolerances are adjusted to satisfy a specified duality gap, but solve the RMP centrally. In the context of dual decomposition in stochastic integer programming, \cite{lubin2013parallelizing} discusses the computational benefits of using tailored interior-point methods such as PIPS-IPM to solve the master problem. Such solvers leverage dual block angular structures by parallelizing the linear algebra.

In privacy-preserving optimization, common techniques include random matrix transformations or perturbations. In \cite{olvi}, the author multiplies the constraints by a randomly generated matrix in the case of equality constraints. A similar method to handle inequality constraints is proposed in \cite{li2013privacy}. In \cite{hong2012secure}, Dantzig-Wolfe decomposition is applied to a transformed problem to preserve privacy, where the constraint matrices of each agent are multiplied by a random matrix, but requires all the data to be centrally located to solve the master problem. Many other methods rely on distributed algorithms to preserve privacy of information \cite{belletti2015privacy, li2020privacy}.

Consensus problems and distributed optimization algorithms have been heavily studied. A classic method is dual decomposition where the linking constraints are relaxed, and the problem becomes separable. The algorithm alternates between solving local subproblems independently and a central step which updates the dual variables as in dual ascent \cite{boyd}. Many variants of dual decomposition have been proposed, such as using subgradients to update the dual iterates when optimizing nonsmooth functions \cite{nedic2009distributed}. Many of the distributed methods are optimized over a network, where agents, treated as nodes, only share limited information with neighbors according to a transition matrix \cite{nedic2010constrained,mota2015distributed,xi2017distributed}. A survey on dual decomposition techniques for distributed optimization and consensus problems is given in \cite{nedic2010survey}. Dual decomposition is known to suffer from weak convergence properties, which led to the augmented Lagrangian method. It provides stronger convergence properties and requires fewer assumptions than dual ascent (the method behind dual decomposition), but leads to an optimization problem that cannot be solved in a distributed way. ADMM has been developed to leverage the decomposability of dual decomposition, and the nice convergence properties of the augmented Lagrangian method \cite{boyd}. Theoretical results and convergence properties of ADMM have been thoroughly studied in \cite{boyd,hong2017linear,nishihara2015general,cai2014direct,beck2017first}. Studies on parameter tuning have also been done, notably the penalty parameter \cite{ghadimi2015optimal,wohlberg2017admm}.

The remainder of the paper is organized as follows. Section \ref{sec:prelims} formally defines the problem structure we are interested in and establishes the notation used throughout. In Section \ref{sec:admm}, we give a brief overview of consensus-based ADMM and discuss traditionally used stopping criteria. In Section \ref{sec:cdwd}, we describe our algorithm and prove bounds on the optimality gap and feasibility violation. We include numerical results to illustrate our method in Section \ref{sec:comp_results}.
\section{Preliminaries} \label{sec:prelims}
We are interested in problems of the form of \eqref{OPT:P}
where $c_n$ and $X_n$ are the cost vector and local constraints of block $n=1,...,N$, and $A_n$ is the constraint matrix of block $n$ in the linking constraints. To simplify notation, we rewrite the master problem as

\begin{mini}
{\lambda}
{\sum_{n=1}^N\sum_{i \in I_n}c_n^i\lambda_{ni}}
{\label{OPT:MP} \tag{$MP$}}
{}
\addConstraint{\sum_{n=1}^N\sum_{i \in I_n}A_n^i\lambda_{ni}}{= t}
\addConstraint{\sum_{i \in I_n} \lambda_{ni}}{= 1,}{\quad n = 1, \dots , N}
\addConstraint{\lambda_{ni}}{\geq 0,}{\quad n = 1, \dots , N}
\end{mini}
where $c_n^i = c_n^\top x_n^i$ and $A_n^i = Ax_n^i$ for all $i \in I_n$ and for all $n$.

We assume each $X_n$ to be a non-empty polytope, so that there exists $L_n>0$ such that $\norm{x_n^i}_2 \leq L_n$ for all extreme points $x_n^i$ of $X_n$. We further assume problem \eqref{OPT:MP} to be feasible and to have an optimal solution. The dual of \eqref{OPT:MP} is
\begin{maxi*}
{\pi, u}
{t^\top \pi + \sum_{n=1}^Nu_n}
{}
{}
\addConstraint{A_n^{i^\top}\pi + u_n}{\leq c_n^i,}{\quad \forall i \in I_n,  \ n = 1, \dots ,N}
\end{maxi*}
where $\pi$ are dual variables associated with the linking constraints and $u_n$ are dual variables associated with the convexity constraints $\sum_{i \in I_n} \lambda_{ni} = 1$ for all $n$. We restrict the presentation of our method and analysis to the case of linking equality constraints without loss of generality.

Let $M > 0$ be an upper bound on the absolute values of the components of $\pi$. We show in Lemma \ref{lem:bound_dual_components} that we can pick $M$ to be finite and polynomial in the entries of the data and get an equivalent dual problem. We will need this constant when proving bounds on optimality gaps. For notational convenience in the proof of Lemma \ref{lem:bound_dual_components}, let $A$ be the horizontal concatenation of the $N$ matrices $A_n$ and let $c$ represent the concatenation of the cost vectors $\{c_n\}_{n=1}^N$. For an index set $J$, let $A_J$ be a matrix formed by columns $j \in J$ of $A$ and $c_J$ be a vector formed by components $j \in J$ of $c$. Let $B \in \mathcal{B}$ be an index set of the columns of $A$ representing a basis, so that $A_B$ is an invertible submatrix of $A$, where $\mathcal{B}$ is the set of all possible bases. Let $\mathcal{N}_B$ be the remaining indices not in $B$. Finally, let $e$ be the vector of all ones of appropriate size and $e_j$ be the vector of all zeros with a one in the $j^{th}$ component.

\begin{lemma}\label{lem:bound_dual_components}
If \eqref{OPT:P} has an optimal solution, then there exists a constant $M > 0$ that is finite and polynomial in the data such that enforcing bound constraints $-M \leq \pi \leq M$ leads to an equivalent dual problem.
\end{lemma}
\begin{proof}
An equivalent big-M formulation of (\ref{OPT:P}) is
\begin{mini}
{x, y^+, y^-}
{\sum_{n=1}^N c_n^\top x_n  + Me^\top y^+ + Me^\top y^-}
{\label{OPT:PM} \tag{$P_M$}}
{}
\addConstraint{\sum_{n=1}^NA_nx_n + y^+ - y^-}{= t}
\addConstraint{x_n}{\in X_n,}{\quad n = 1, \dots, N}
\addConstraint{y^+, y^-}{\geq 0.}
\end{mini}
Since \eqref{OPT:P} is assumed to have an optimal solution, \eqref{OPT:PM} is equivalent to \eqref{OPT:P} if $M$ is large enough and at any optimal solution, $y^+ = y^- = 0$. Note that the dual of the master problem reformulated from \eqref{OPT:PM} is the same as the dual of \eqref{OPT:MP} with added bounds $-M \leq \pi \leq M$.

If we assume problem \eqref{OPT:PM} is solved using the simplex method where at any point in the algorithm, only columns with a negative reduced cost can enter the basis, we require $M$ to be big enough so that reduced costs of $y^+$ and $y^-$ remain positive for all bases $B \in \mathcal{B}$ (this is not necessary but sufficient for the proof of existence). At a basis $B$ which does not contain a component from either $y^+$ or $y^-$, the reduced costs of components $y^+_j$ and $y^-_j$ are $M - c_B^\top A_B^{-1}e_j$ and $M + c_B^\top A_B^{-1}e_j$, respectively, for $j \in \mathcal{N}_B$.

We need $M - c_B^\top A_B^{-1}e_j > 0$ and $M + c_B^\top A_B^{-1}e_j > 0$, or simply $M >\abs{c_B^\top A_B^{-1}e_j}$. 
Thus, $M >\abs{c_B^\top A_B^{-1}e_j}$ $\forall B \in \mathcal{B}$ and $\forall j \in \mathcal{N}_B$ would be sufficient to obtain an equivalent problem and the result follows.
\end{proof}


We refer to the optimal objective values of the master problem and its dual by $z_{MP}^*$ and $z_{DM}^*$. When considering a subset of the extreme points, we refer to the restricted primal and dual problems as RMP and RDM, and their optimal values by $z_{RMP}^*$ and $z_{RDM}^*$, respectively. Approximate solutions and their objective values are denoted by a hat. The notations $\norm{\cdot}$ and $\norm{\cdot}_F$ refer to the $\ell_2$-norm and Frobenius norm, respectively. Finally, the terms agent and block will be used interchangeably.

\section{ADMM Overview}\label{sec:admm}
We give a brief overview of consensus-based ADMM. We present the algorithm and known convergence conditions. Detailed discussion of the ADMM method and its convergence properties can be found in \cite{boyd}.

\subsection{Consensus-Based ADMM}
Consensus-based ADMM is well-suited for problems of the form
\begin{maxi}
{x}
{\sum_{n=1}^Nf_n(x)}
{\label{OPT:admm} \tag{A1}}
{}
\addConstraint{A_nx}{\leq b_n}{\quad \forall n=1, \dots ,N}
\end{maxi}
where $f_n : \mathbb{R}^n \rightarrow \mathbb{R}$ are convex, proper and closed functions, and $A \in \mathbb{R}^{m \times n}$. The objective function and constraints are linked solely by the variable $x$. We can equivalently rewrite \eqref{OPT:admm} as
\begin{maxi!}
{x_n, x}
{\sum_{n=1}^Nf_n(x_n)  \nonumber}
{\label{OPT:admm_sep}}
{}
\addConstraint{A_nx_n}{\leq b_n}{\quad n = 1, \dots ,N} \label{A_const}
\addConstraint{x_n}{= x}{\quad n = 1, \dots ,N.} \label{copy_const}
\end{maxi!}

Let $\lambda_n \in \mathbb{R}^m$ and $\alpha_n \in \mathbb{R}^n$ be the Lagrangian multipliers of \eqref{A_const} and \eqref{copy_const}, respectively, for $n = 1, \dots ,N$. Taking the augmented Lagrangian of \eqref{OPT:admm_sep} gives
\begin{maxi}
{x_n, x}
{\sum_{n=1}^N \left[ f_n(x_n) + \alpha_n^\top (x - x_n) - \frac{\rho}{2}\norm{x - x_n}^2 \right]}
{\label{OPT:admm_aug} \tag{AL}}
{}
\addConstraint{A_n x_n}{\leq b_n}{\quad n = 1, \dots ,N}
\end{maxi}
where $\rho > 0$ is a predetermined penalty parameter. Define the objective of  \eqref{OPT:admm_aug} as
\begin{equation*}
\mathcal{L}_\rho(x,x_1,...,x_N,\alpha_1,...,\alpha_N) = \sum_{n=1}^N \left[ f_n(x_n) + \alpha_n^\top (x - x_n) - \frac{\rho}{2}\norm{x - x_n}^2 \right]
\end{equation*}
The ADMM method consists of alternating between maximizing the function $\mathcal{L}_\rho$ over $(x,x_1,...,x_N)$ and minimizing over $(\alpha_1,...,\alpha_N)$, where the maximization step is done sequentially, so that we first maximize over $(x_1,...x_N)$ before maximizing over $x$. This allows to solve the former in a distributed fashion. The ADMM steps at an iteration $k$ can be summarized as follows:
\begin{subequations} \label{eq:admm_steps}
\begin{align}
&x_n^{k+1} \leftarrow \underset{x_n}{\text{argmax }}\mathcal{L}_\rho(x^k,x_1,...,x_N,\alpha_1^k,...,\alpha_N^k), \ \forall n = 1,...,N  \\
&x^{k+1} \leftarrow \underset{x}{\text{argmax }}\mathcal{L}_\rho(x,x_1^{k+1},...,x_N^{k+1},\alpha_1^k,...,\alpha_N^k) \label{subeq:admm_step2} \\
&\alpha_n^{k+1} \leftarrow \alpha_n^k - \rho(x^{k+1} - x_n^{k+1}), \ \forall n = 1,...,N  \label{subeq:admm_step3}
\end{align}
\end{subequations}
Note that \eqref{subeq:admm_step3} is a gradient step where the step size is the penalty parameter $\rho$ and \eqref{subeq:admm_step2} is an unconstrained maximization problem for which there exists a closed form solution. We have:
\begin{align*}
&\nabla_x \mathcal{L}_\rho(x,x_1^{k+1},...,x_N^{k+1},\alpha_1^k,...,\alpha_N^k) = 0 \\
\Rightarrow &\sum_{n=1}^N \left[ \alpha_n^k - \rho(x^{k+1} - x_n^{k+1}) \right] = 0 \\
\Rightarrow & x^{k+1} = \frac{1}{N}\sum_{n=1}^Nx_n^{k+1} + \frac{1}{N\rho}\sum_{n=1}^N\alpha_n^k.
\end{align*}

\subsection{Convergence and Stopping Criteria}
We first note that dual feasibility in \eqref{OPT:admm_sep} is equivalent to $\nabla_{x_n} f_n(x_n) - A^\top_n\lambda_n - \alpha_n = 0$. From the optimality conditions of \eqref{OPT:admm_aug}, we have at iteration $k$:
\begin{align*}
&\nabla_{x_n} f_n(x_n^{k+1}) - A^\top_n\lambda_n^{k+1} - \alpha_n^k + \rho(x^k-x_n^{k+1}) = 0 \\
\Rightarrow &\nabla_{x_n} f_n(x_n^{k+1}) - A^\top_n\lambda_n^{k+1} - \alpha_n^k + \rho(x^k-x_n^{k+1} +x^{k+1} - x^{k+1}) = 0 \\
\Rightarrow &\nabla_{x_n} f_n(x_n^{k+1}) - A^\top_n\lambda_n^{k+1} - \alpha_n^k + \rho(x^{k+1}-x_n^{k+1}) + \rho(x^{k} - x^{k+1}) = 0 \\
\Rightarrow &\nabla_{x_n} f_n(x_n^{k+1}) - A^\top_n\lambda_n^{k+1} - \alpha_n^{k+1} = \rho(x^{k+1} - x^{k}).
\end{align*}
Thus, dual feasibility in \eqref{OPT:admm_sep} amounts to having $\rho(x^{k+1} - x^{k}) = 0$.

Let $\boldsymbol{\alpha}$ be the vertical concatenation of vectors $\{\alpha_n\}_{n=1}^N$. As proven in \cite{boyd}, under the assumption that the functions $f_n$ in \eqref{OPT:admm} are convex, proper and closed, and assuming that $\mathcal{L}_0(x,x_1,...,x_N,\boldsymbol{\alpha})$, where $\rho = 0$, has a saddle point, then as $k \rightarrow \infty$, we have the following:
\begin{enumerate}[(i)]
\item Primal feasibility violation vanishes: $\sqrt{\sum_{n=1}^N\norm{x^{k+1} - x^{k+1}_n}^2} \rightarrow 0, \  n=1,...,N$
\item Dual feasibility violation vanishes: $\rho \norm{x^{k+1} - x^{k}} \rightarrow 0$
\item Optimality gap vanishes: $\norm{f(x^*) - f(x^k)} \rightarrow 0$
\item Dual vector $\boldsymbol{\alpha}$ converges to an optimal dual solution: $\norm{\boldsymbol{\alpha^k} - \boldsymbol{\alpha}^*} \rightarrow 0$
\end{enumerate}
In (i), we define primal feasibility violation to be
\begin{align*}
\begin{Vmatrix} x^{k+1} - x^{k+1}_1 \\ \vdots \\ x^{k+1} - x^{k+1}_N\end{Vmatrix} = \sqrt{\sum_{n=1}^N\norm{x^{k+1} - x^{k+1}_n}^2}
\end{align*}
This implies that we reach consensus as $k \rightarrow \infty$, i.e $\norm{x^{k+1} - x^{k+1}_n} \rightarrow 0$ for all $n$.

For our purposes, we also assume the functions $f_n$ to be differentiable. This assumption is satisfied in our case since we are dealing with linear cost functions.  As suggested in \cite{boyd}, it is reasonable to terminate ADMM once we reach primal and dual feasibility within some tolerance. Given specified tolerances $\epsilon_p$ and $\epsilon_d$, we terminate ADMM once $\sqrt{\sum_{n=1}^N\norm{x^{k+1} - x^{k+1}_n}^2} \leq \epsilon_p$ and $\rho\norm{x^{k+1} - x^{k}} \leq \epsilon_d$.

\section{Consensus-Based Dantzig-Wolfe Algorithm}\label{sec:cdwd}
We first present the consensus-based Dantzig-Wolfe decomposition (CDWD) algorithm before deriving error bounds on the optimality gap and feasibility violation.

\subsection{CDWD Algorithm}\label{subsec:cdwd_alg}
We define $k$ to be the ADMM iteration counter and $\ell$ to be the Dantzig-Wolfe outer iteration counter. To solve the restricted master problem (RMP) in a distributed fashion, we  solve a reformulation of the dual of RMP. The reformulation permits us to perform consensus-based ADMM. We split the dual vector $\pi$ associated with the linking constraints into $N$ copies as in \eqref{OPT:admm_sep} to get the following equivalent formulation:

\begin{maxi}
{\pi_n, u_n}
{\sum_{n=1}^N \left[ \frac{1}{N} t^\top \pi_n + u_n \right]}
{\label{OPT:DM}\tag{$DM$}}
{}
\addConstraint{A_n^{i^\top}\pi_n + u_n}{\leq c_n^i}{\quad \forall i \in I_n, \ n = 1, \dots ,N}
\addConstraint{\pi_n}{= \pi}{\quad n = 1,...,N.}
\end{maxi}
If there exists linking inequality constraints, then we can simply add appropriate non-negativity or non-positivity constraints on $\pi_n$ for all $n$, and the remaining steps of the algorithm follow.

Note that in \eqref{OPT:MP}, the problem is linked by the rows. Performing ADMM directly on \eqref{OPT:MP} would lead to $N$ blocks where we would need to optimize with respect to each $x_n$ sequentially. Not only is the problem no longer decomposable, but additional conditions and correction steps are necessary to ensure convergence of ADMM with more than two blocks \cite{chen2016direct, deng2017parallel, goncalves2019two}. The dual of the master problem, however, is linked by the decision variable $\pi$. Performing ADMM on \eqref{OPT:DM} leads to a more natural consensus-based algorithm with a guaranteed convergence, where the first ADMM block corresponds to solving for each $\pi_n$ independently, and the second block corresponds to optimizing with respect to $\pi$. As a result, working on \eqref{OPT:DM} avoids any potential violation of privacy, as agents need only share dual vectors with a central coordinator. Although one could solve \eqref{OPT:MP} using a dual decomposition method with subgradients while maintaining decomposability, two problems arise. First, privacy could be violated when computing the subgradients (i.e. the residuals of the linking constraints) \cite{lou2017privacy, yan2012distributed}. Second, using an averaging scheme would be necessary to recover an approximate primal solution. Indeed, it is well-known that simply solving the Lagrangian dual at the optimal dual solution does not necessarily return a feasible primal solution \cite{frangioni2005lagrangian, boyd}. Although out of the scope of this paper and further experiments are required, we note that averaging schemes did not perform well for simple instances in our setting. We show that by performing ADMM on \eqref{OPT:DM}, we circumvent primal-recovery issues and can easily retrieve high quality primal solutions via the Lagrangian multipliers associated with the constraints in \eqref{OPT:ARDMn} (defined below).

We denote the restricted problem of \eqref{OPT:DM}, i.e. one involving constraints corresponding to only a subset of the columns, by RDM and its optimal value by $z_{RDM}^*$. We take the augmented Lagrangian of RDM by relaxing the copy constraints as in (\ref{OPT:admm_aug}) and get a separable problem with respect to variables $(\pi_n,u_n)$:
\begin{maxi*}
{\pi_n, u_n}
{\sum_{n=1}^N \left[\frac{1}{N} t^\top\pi_n + u_n + \alpha_n^\top(\pi - \pi_n) - \frac{\rho}{2}\norm{\pi - \pi_n}^2 \right]}
{}
{}
\addConstraint{A_n^{i^\top}\pi_n + u_n}{\leq c^i_n}{\quad \forall i \in I_n, \ n = 1, \dots ,N.}
\end{maxi*}

At iteration $k$ of ADMM and using current iterates $\pi^k$ and $\alpha_n^k$, each agent $n$ solves
\begin{align}\label{OPT:ARDMn}
\tag{\text{$ARDM_n$}}
\begin{split}
\max_{\pi_n,u_n} &\qquad \frac{1}{N} t^\top\pi_n + u_n + \alpha_n^{k^\top}(\pi^k - \pi_n) - \frac{\rho}{2}\norm{\pi^k - \pi_n}^2 \\
\text{s.t} &\qquad A_n^{i^\top}\pi_n + u_n \leq c^i_n, \ \forall i \in I^\ell_n 
\end{split}
\end{align}
where $I_n^\ell \subseteq I_n$ is the index set of extreme points of block $n$ at outer iteration $\ell$. From \eqref{eq:admm_steps}, the steps to solving RDM can be summarized as follows:
\begin{enumerate}[1.]
\item Each agent solves \eqref{OPT:ARDMn} and collects optimal solutions $(\pi_n^{k+1},u^{k+1}_n)$
\item $\pi^{k+1} \leftarrow \frac{1}{N}\sum_{n=1}^N(\pi_n^{k+1}) + \frac{1}{N\rho}		\sum_{n=1}^N\alpha_n^k$
\item $\alpha^{k+1}_n = \alpha^k_n - \rho(\pi^{k+1} - \pi^{k+1}_n)$
\end{enumerate}

First note that $A_n^{i^\top}\pi_n^{k+1} + u_n^{k+1} \leq c_n^i$ is satisfied for all $i \in I_n^\ell$ and for all $n$, since $(\pi_n^{k+1},u_n^{k+1})$ is a solution of \eqref{OPT:ARDMn}. Thus, $\pi_n^{k+1} = \pi^{k+1}$ are the only violated constraints. To avoid confusion, we refer to $\sqrt{\sum_{n=1}^N\norm{\pi^{k+1} - \pi^{k+1}_n}^2}$ as the dual feasibility violation, and $\rho\norm{\pi^{k+1} - \pi^k}$ as the primal feasibility violation. Note that this is the opposite of what is defined in Section \ref{sec:admm} because we are performing ADMM on the dual problem here. We then perform steps 1-3 until $\sqrt{\sum_{n=1}^N\norm{\pi^{k+1} - \pi^{k+1}_n}^2} \leq \epsilon_d$ and $\rho\norm{\pi^{k+1} - \pi^k} \leq \epsilon_p$, where $\epsilon_d$ and $\epsilon_p$ are dual and primal feasibility tolerances, respectively. Each agent  $n$ then solves a pricing subproblem to look for an extreme point with negative reduced cost:
\begin{align*}
z_{SEP}^n = \min_x\{c^\top_nx_n - \pi^{{k+1}^\top}A_nx_n - u^{k+1}_n: x_n \in X_n\}.
\end{align*}

Let $x^*_n$ be an optimal solution. In standard DWD, we would add $x^*_n$ as a new column if $z_{SEP}^n < 0$. However, the dual solution $(\pi^{k+1},\{u_n^{k+1}\}_{n=1}^N)$ is $\epsilon$-optimal and only close to feasible for the current RMP. It is possible that we find a column whose reduced cost is negative and close to 0 when evaluated at the approximate dual solutions, but is in fact already in the current RMP. It is also possible that at the (unavailable) optimal dual solution, the reduced cost is actually positive and the extreme point should not be added. To ensure a finite algorithm, agent $n$ only adds $x^*_n$ as a new extreme point if $z_{SEP}^n < -\max_{i \in I_n^\ell}\{\norm{A_n^i}\}\epsilon_d$ . In Lemma \ref{lem:sep_bounds1}, we show that ADMM terminates with $c^i_n - A_n^{i^\top}\pi^{k+1} - u_n^{k+1} \geq -\norm{A_n^i}\epsilon_d$ for all $i$ and $n$. Thus, if $-\max_{i \in I_n^\ell}\{\norm{A_n^i}\}\epsilon_d \leq z^n_{SEP} = c_n^\top x_n^* - \pi^{{k+1}^\top}A_nx_n^* - u_n^{k+1} < 0$, then we cannot guarantee that $x_n^*$ is a necessary extreme point. In other words, we can only trust $z^n_{SEP}$ within $\max_{i \in I_n^\ell}\{\norm{A_n^i}\}\epsilon_d$. This is necessary for the analysis in Section \ref{subsec:convergence}.

\begin{lemma}\label{lem:sep_bounds1}
At outer-iteration $\ell$, if ADMM terminates with $\norm{\pi^{k+1} - \pi_n^{k+1}} \leq \epsilon_d$ for all $n$, we have $$c^i_n - A_n^{i^\top}\pi^{k+1} - u_n^{k+1} \geq -\norm{A_n^i}\epsilon_d$$ for all $i \in I_n^\ell, n = 1,...,N$.
\end{lemma}
\begin{proof}
We have
\begin{alignat*}{3}
& &&\begin{Vmatrix} \pi^{k+1} - \pi^{k+1}_1 \\ \vdots \\ \pi^{k+1} - \pi^{k+1}_N\end{Vmatrix} && \leq \epsilon_d \\
&\Rightarrow \sum_{n=1}^N \ &&\norm{\pi^{k+1} - \pi^{k+1}_n}^2 &&\leq \epsilon_d^2 \\
&\Rightarrow &&\norm{\pi^{k+1} - \pi^{k+1}_n}^2 &&\leq \epsilon_d^2 \\
&\Rightarrow &&\norm{\pi^{k+1} - \pi^{k+1}_n} &&\leq \epsilon_d, \ \forall n = 1, \dots ,N.
\end{alignat*}
For any $n$, computing the distance between $c_n^i - A_n^{i^\top}\pi^{k+1} - u_n^{k+1} $ and $c_n^i - A_n^{i^\top}\pi^{k+1}_n - u_n^{k+1} $ gives us
\begin{align}\label{rcost_bounds}
\norm{c_n^i - A_n^{i^\top}\pi^{k+1} - u_n^{k+1} - c_n^i + A_n^{i^\top}\pi^{k+1}_n + u_n^{k+1}} &= \norm{A_n^{i^\top}(\pi^{k+1} - \pi^{k+1}_n)} \\
&\leq \norm{A_n^i}\epsilon_d, \ \forall i \in I_n^\ell\nonumber
\end{align}
Since $c_n^i - A_n^{i^\top}\pi^{k+1}_n - u_n^{k+1}\geq 0$ for all $i \in I_n^\ell$, (\ref{rcost_bounds}) implies $c^i_n - A_n^{i^\top}\pi^{k+1} - u_n^{k+1} \geq -\norm{A_n^i}\epsilon_d$ for all $i \in I_n^\ell$ and $n$.
\end{proof}

Once the columns are added to the RMP, we use the solutions of the last iterates $\pi^{k+1}$ and $\alpha_n^{k+1}$ as warm starts for $\pi^1$ and $\alpha_n^1$ for all $n$ in the next outer iteration $\ell + 1$. If $z_{SEP}^n \geq -\max_{i \in I_n^{\ell+1}}\{\norm{A_n^i}\}\epsilon_d$ for all $n$, we terminate the algorithm and retrieve approximate primal solutions $\hat{x}_n \leftarrow \sum_{i \in I^{\ell+1}_n}\lambda_{ni}^{k+1}x^i_n$ for all  $n = 1,...,N$, where $\lambda_{ni}^{k+1}$ are the Lagrangian multipliers associated with the constraints in \ref{OPT:ARDMn}, $i \in I_n^{\ell+1}, n = 1,...,N$. The CDWD algorithm is summarized in Algorithm \ref{alg:cdwd}. We describe the algorithm in a master-worker framework where agents are referred to as processors. At each ADMM iteration, the master node calls the \texttt{BROADCAST()} function to send the current estimate of $\pi$ to each processor, and the \texttt{RECEIVE()} function to collect each processor's dual solution $\pi_n$ obtained from solving \ref{OPT:ARDMn}.

Note that for computational benefits, we repeatedly solve Algorithm \ref{alg:cdwd} where we start with loose tolerances and dynamically tighten them to achieve a desired accuracy, but omit this discussion here for ease of exposition as it does not affect the analysis. 

\begin{algorithm}
\caption{CDWD Algorithm}
\label{alg:cdwd}
\begin{algorithmic}[1]
\State Input: tolerances $\epsilon_p, \epsilon_d \geq 0$, penalty parameter $\rho > 0$
\State Let $I^1_n$ be the initial set of columns for each block $n$ 
\State Initialize $\pi^1,\alpha_n^1$ for all $n$ and $\ell = 0$
\While{$I_n^{\ell+1} \neq I_n^\ell$ for some $n$}
	\State $\ell \leftarrow \ell + 1$
	\State Initialize primal and dual residuals $r_p = \infty$ and $r_d = \infty$
	\item[]
	\State /*Solve RDM using consensus-based ADMM*/
	\State Initialize $k=0$
	\While{$r_d > \epsilon_d$ and $r_p > \epsilon_p$}
		\State $k \leftarrow k+1$
		\State \texttt{BROADCAST($\pi^{k}$)}	
		\ForEach{\text{processor} $n = 1,...,N$}
			\State Solve (\ref{OPT:ARDMn})
			\State Collect optimal solutions $(\pi_n^{k+1},u^{k+1}_n)$
			\State Collect Lagrangian multipliers $\lambda^{k+1}_n$
		\EndFor
		\State \texttt{RECEIVE($\{\pi^{k+1}_n\}_{n = 1}^N$)}
		\State $\pi^{k+1} \leftarrow \frac{1}{N}\sum_{n=1}^N(\pi_n^{k+1}) + \frac{1}{N\rho}		\sum_{n=1}^N\alpha_n^k$
		\State $\alpha^{k+1}_n = \alpha^k_n - \rho(\pi^{k+1} - \pi^{k+1}_n)$
		\State $r_d \leftarrow \sqrt{\sum_{n=1}^N\norm{\pi^{k+1}-\pi_n^{k+1}}^2}$
		\State $r_p \leftarrow \rho\norm{\pi^{k+1}-\pi^k}$
	\EndWhile
	\State \texttt{BROADCAST($\pi^{k+1}$)}
	\item[]
	\State /*Solve pricing subproblems*/
	\ForEach{\text{processor} $n = 1,...,N$}
		\State $z_{SEP}^n \leftarrow \min_{x_n}\{c^\top_nx_n - \pi^{{k+1}^\top}A_nx_n - u^{k+1}_n: x_n \in X_n\}$
		\State Let $x_n^i$ be the optimal solution
		\If{$z_{SEP}^n < -\max_{i \in I_n^\ell}\{\norm{A_n^i}\} \epsilon_d$}
			\State Add extreme point $x^i_n$: $I_n^{\ell+1} \leftarrow I_n^\ell \cup \{i\}$
		\Else
			\State $I_n^{\ell+1} \leftarrow I_n^\ell$
		\EndIf
	\EndFor
	\State $\pi^1 \leftarrow \pi^{k+1}$
	\State $\alpha_n^1 \leftarrow \alpha_n^{k+1}, \ \forall n = 1,...,N$
\EndWhile
\item[]
\State /* Each processor $n$ retrieves primal solution*/
\State $\hat{x}_n \leftarrow \sum_{i \in I^{\ell+1}_n}\lambda_{ni}^{k+1}x^i_n, \ \forall n = 1,...,N$
\end{algorithmic}
\end{algorithm}

\subsection{Convergence}\label{subsec:convergence}
We now prove the convergence of CDWD and provide bounds on the optimality gap and feasibility violation. The quality of the dual solutions obtained by the consensus ADMM algorithm directly affects the quality of the recovered primal solution. We are able to reduce the optimality gap and feasibility violation by tweaking the primal and dual infeasibility tolerances $\epsilon_p$ and $\epsilon_d$. Recall since we are solving the dual of (\ref{OPT:MP}) using ADMM, we refer to the Lagrangian multipliers $\alpha_n$ in the objective of (\ref{OPT:ARDMn}) and the multipliers $\lambda_{ni}$ associated with the constraints as primal variables, and $\pi, \pi_n$ and $u_n$ as dual variables; we refer to $\sqrt{\sum_{n=1}^N \norm{\pi^{k+1} - \pi^{k+1}_n}^2}$ as the dual feasibility violation and $\rho\norm{\pi^{k+1} - \pi^k}$ as the primal feasibility violation. Moreover, recall that $z_{MP}^*$ and $z_{DM}^*$ refer to the optimal objective values of the master problem (\ref{OPT:MP}) and its dual, respectively; $z_{RMP}^*$ and $z_{RDM}^*$ refer to the optimal values of their restrictive counterparts; objective values and solutions resulting from the CDWD algorithm are denoted by a hat such as $\hat{z}_{RDM}$. 

As shown in \cite{boyd} and other sources in the literature, given tolerances $\epsilon, \epsilon_p, \epsilon_d > 0$, we can assume that ADMM terminates with $z_{RDM}^* - \hat{z}_{RDM} \leq \epsilon$, $\rho\norm{\pi^{k+1} - \pi^k} \leq \epsilon_p$ and $\sqrt{\sum_{n=1}^N \norm{\pi^{k+1} - \pi^{k+1}_n}^2} \leq \epsilon_d$.
The following lemmas will be helpful in proving the error bounds.

\begin{lemma} \label{lem:sum_alphas}
After the first iteration of CDWD, the Lagrangian multipliers    $\alpha^k_n$  associated with the copy constraints are primal feasible for all $n$, i.e for $k \geq 0$, we have $\sum_{n=1}^N \alpha_n^{k+1} = 0$.
\end{lemma}

\begin{proof}[Proof]
From the updates, we have with $k \geq 0$:
\begin{align*}
&\pi^{k+1} = \frac{1}{N}\sum_{n=1}^N \pi_n^{k+1} + \frac{1}{N\rho}\sum_{n=1}^N\alpha_n^k \\
\Rightarrow &\sum_{n=1}^N\alpha^k_n = N\rho\pi^{k+1} - \rho \sum_{n=1}^N\pi_n^{k+1} \\
\Rightarrow &\sum_{n=1}^N\alpha^k_n = \rho \left(N\pi^{k+1} - \sum_{n=1}^N \pi_n^{k+1} \right)
\end{align*}
and
\begin{align*}
&\alpha^{k+1}_n = \alpha^k_n - \rho(\pi^{k+1} - \pi^{k+1}_n), \ \forall n = 1,...,N.
\end{align*}
Summing over $n$, we get
\begin{align*}
\sum_{n=1}^N\alpha^{k+1}_n &= \sum_{n=1}^N\alpha^{k}_n - \rho \left( N\pi^{k+1} - \sum_{n=1}^N\pi_n^{k+1} \right) \\
&=\rho \left(N\pi^{k+1} - \sum_{n=1}^N \pi_n^{k+1} \right) - \rho \left( N\pi^{k+1} - \sum_{n=1}^N\pi_n^{k+1}\right) \\
&= 0.
\end{align*}
\end{proof}

Theorem \ref{thm1} establishes the feasibility violation at the recovered primal solution.
\begin{theorem}[Feasibility Violation]\label{thm1}
Given a primal feasibility tolerance $\epsilon_p > 0$, CDWD terminates with a solution $\hat{x}_n = \sum_{i \in I_n^{\ell+1}}\lambda_{ni}^{k+1}x_n^i$ such that:
\begin{align*}
&\norm{\sum_{n = 1}^NA_n\hat{x}_n - t} \leq N\epsilon_p \\
&\sum_{i \in I_n^\ell}\lambda^{k+1}_{ni} = 1, \ \forall n
\end{align*}
\end{theorem}

\begin{proof}
At outer iteration $\ell$ and iteration $k$ of ADMM, let the Lagrangian functions of (\ref{OPT:ARDMn}) for each $n$ be:
\begin{multline*}
Q_n(\pi_n, \alpha_n^k,\pi^k,\lambda) = \frac{1}{N}t^\top\pi_n + u_n + \alpha_n^{k^\top}(\pi^k - \pi_n) - \frac{\rho}{2} \norm{\pi^k - \pi_n}_2^2
	+ \sum_{i \in I_n^\ell} \lambda_{ni}(c_n^i - A_n^{i^\top} \pi_n - u_n)
\end{multline*}
where $\{\lambda_{ni}\}_{i \in I^\ell_n}$ are the multipliers of the constraints in (\ref{OPT:ARDMn}).

We have the following optimality conditions in (\ref{OPT:ARDMn}):
\begin{align*}
&\lambda^{k+1}_{ni}(c_n^i - A_n^{i^\top} \pi^{k+1}_n - u^{k+1}_n) = 0,\ \forall i\in I_n^\ell \tag{Complementary Slackness} \\
&\lambda^{k+1}_{ni} \geq 0, \ \forall i\in I_n^\ell \tag{Dual Feasibility} \\
\tag{Stationarity}
\begin{split}
&\nabla_{\pi_n}Q_n = \frac{1}{N}t - \alpha^k_n + \rho(\pi^k - \pi^{k+1}_n) - \sum_{i \in I_n^\ell}A_n^i \lambda^{k+1}_{ni} = 0  \\
&\nabla_{u_n}Q_n = 1 - \sum_{i \in I_n^\ell}\lambda^{k+1}_{ni} = 0.
\end{split}
\end{align*}

Thus, the convexity constraints in RMP $\sum_{i \in I_n^\ell}\lambda^{k+1}_{ni} = 1$ are satisfied for all $n$ and $\lambda^{k+1}_{ni} \geq 0$ for all $n$ and $i$.

We rewrite the stationarity condition with respect to $\pi_n$ as
\begin{align}
\nabla_{\pi_n}Q_n &= \frac{1}{N}t - \alpha^k_n + \rho(\pi^k - \pi^{k+1}_n) - \sum_{i \in I_n^\ell}A_n^i \lambda^{k+1}_{ni} \nonumber\\
&= \frac{1}{N}t - \alpha^k_n + \rho(\pi^k - \pi^{k+1} + \pi^{k+1} - \pi^{k+1}_n) - \sum_{i \in I_n^\ell}A_n^i \lambda^{k+1}_{ni} \nonumber\\
&= \frac{1}{N}t - \alpha^k_n + \rho(\pi^k - \pi^{k+1}) + \rho(\pi^{k+1} - \pi^{k+1}_n) - \sum_{i \in I_n^\ell}A_n^i \lambda^{k+1}_{ni} \nonumber\\
&= \frac{1}{N}t - \rho(\pi^{k+1} - \pi^k) - \sum_{i \in I_n^\ell}A_n^i \lambda_{ni}^{k+1} - \alpha_n^{k+1} \label{feas_bound_split}
\end{align}
where the last equality holds from $\alpha^{k+1}_n = \alpha^k_n - \rho(\pi^{k+1} - \pi^{k+1}_n)$.

Summing $\nabla_{\pi_n}Q_n$ over $n = 1,...,N$, we get
\begin{align*}
\sum_{n = 1}^N\nabla_{\pi_n}Q_n &= t  - \rho \sum_{n=1}^N(\pi^{k+1} - \pi^k) - \sum_{n = 1}^N\sum_{i \in I_n^\ell}A_n^i\lambda^{k+1}_{ni} - \sum_{n = 1}^N \alpha^{k+1}_n\\
&= t - \rho N(\pi^{k+1}-\pi^k) - \sum_{n = 1}^N\sum_{i \in I_n^\ell}A_n^i\lambda^{k+1}_{ni} \\
&= 0
\end{align*}
where the second equality follows because $\sum_{n = 1}^N \alpha^{k+1}_n = 0$ from \ref{lem:sum_alphas}.

Then
\begin{align*}
&\norm{\sum_{n = 1}^N\nabla_{\pi_n}Q_n} = 0 \geq \norm{t - \sum_{n = 1}^N\sum_{i \in I_n^\ell}A_n^i\lambda^{k+1}_{ni}} - N\rho\norm{\pi^{k+1} - \pi^k} \\
\Rightarrow &\norm{t - \sum_{n = 1}^N\sum_{i \in I_n^\ell}A_n^i\lambda^{k+1}_{ni}} \leq N\rho\norm{\pi^{k+1} - \pi^k} \\
\Rightarrow &\norm{t - \sum_{n = 1}^NA_n\hat{x}_n} \leq N\epsilon_p
\end{align*}
where $\hat{x}_n = \sum_{i \in I_n^\ell}\lambda_{ni}^{k+1}x_n^i = \sum_{i \in I_n^{\ell+1}}\lambda_{ni}^{k+1}x_n^i$, since $I_n^{\ell+1} = I_n^\ell$ for all $n$ when CDWD terminates.
\end{proof}

Before deriving error bounds on the optimality gap, we first introduce bounds on $z^*_{DM} - \hat{z}_{RDM}$ and $\hat{z}_{RMP} - \hat{z}_{RDM}$. Adding these two will then give us bounds on $\hat{z}_{RMP} - z^*_{DM}$, or equivalently $\hat{z}_{RMP} - z^*_{MP}$. The following is a known relationship between $\hat{z}_{RDM}$, $z^*_{DM}$ and $z^*_{RDM}$ (cf.~\cite{lubbecke2005}).
\begin{lemma}\label{lem:Rdual_obj_bounds}
After terminating ADMM, we have $\hat{z}_{RDM} + \sum_{n = 1}^N \min\{0,z^n_{SEP}\} \leq z^*_{DM} \leq z^*_{RDM}$.
\end{lemma}
\begin{proof}
If $z_{SEP}^n < 0$ for some $n$, then we can set $\hat{u}_n' = \hat{u}_n + z_{SEP}^n$. Doing so for each $n$, we get a feasible solution $(\hat{\pi},\{\hat{u}_n'\}_{n=1}^N)$ to (\ref{OPT:DM}) with objective value $\hat{z}_{RDM} + \sum_{n = 1}^N \min\{0,z^n_{SEP}\} \leq z^*_{DM}$. Moreover, $z_{DM}^* \leq z_{RDM}^*$ since RDM is a relaxation of (\ref{OPT:DM}).
\end{proof}

\begin{proposition}\label{prop:gap_dual}
Given that ADMM terminates with $z_{RDM}^* - \hat{z}_{RDM} \leq \epsilon$ and we terminate CDWD when $z_{SEP}^n \geq -\max_{i \in I_n^\ell}\{\norm{A_n^i}\}\epsilon_d$ for all $n$, we have
\begin{align*}
-\epsilon_d\sum_{n=1}^N \norm{A_n}_FL_n \leq z_{DM}^* - \hat{z}_{RDM} \leq \epsilon
\end{align*} 
where $L_n$ is a bound on all extreme points of the set $X_n$, defining the local constraints of block $n$.
\end{proposition}

\begin{proof}
By Lemma \ref{lem:Rdual_obj_bounds}, $\hat{z}_{RDM} + \sum_{n = 1}^N \min\{0,z^n_{SEP}\} \leq z^*_{DM} \leq z^*_{RDM}$. Since $z^*_{RDM} - \hat{z}_{RDM} \leq \epsilon$ and $z_{SEP}^n \geq -\max_{i \in I_n^\ell}\{\norm{A_n^i}\}\epsilon_d$ for all $n$ after terminating ADMM, we have:
\begin{align*}
&\hat{z}_{RDM} + \sum_{n = 1}^N \min\{0,z^n_{SEP}\} \leq z^*_{DM} \leq z^*_{RDM} \\
\Rightarrow &\hat{z}_{RDM} - \hat{z}_{RDM} + \sum_{n = 1}^N \min\{0,z^n_{SEP}\} \leq z^*_{DM} - \hat{z}_{RDM} \leq z^*_{RDM} - \hat{z}_{RDM} \\
\Rightarrow &-\sum_{n=1}^N\max_{i \in I_n^\ell}\{\norm{A_n^i}\}\epsilon_d \leq z^*_{DM} - \hat{z}_{RDM} \leq \epsilon \\
\Rightarrow &-\epsilon_d\sum_{n=1}^N \norm{A_n}_FL_n \leq z^*_{DM} - \hat{z}_{RDM} \leq \epsilon \end{align*}
where the last set of inequalities holds from our assumption that $\norm{x_n^i} \leq L_n$ for all extreme points of block $n$: $\norm{A_n^i} = \norm{A_nx_n^i} \leq \norm{A_n}_FL_n$.
\end{proof}

\begin{proposition}\label{prop:restricted_duality_gap}
Terminating ADMM with primal and dual feasibility tolerances $\epsilon_p$ and $\epsilon_d$, respectively, we have at any outer iteration $\ell$
\begin{align*}
\abs{\hat{z}_{RMP} - \hat{z}_{RDM}} \leq \epsilon_d \sum_{n=1}^N \norm{A_n}_F L_n + mM N\epsilon_p
\end{align*}
where $m$ is the number of linking constraints, $M$ is an upperbound on the absolute values of the components of $\pi$ as in Lemma \ref{lem:bound_dual_components}, and $N$ is the number of blocks.
\end{proposition}
\begin{proof}
The complementary slackness conditions for RDM are
\begin{align}
&\lambda_{ni}(c_n^i-A_n^{i^\top}\pi - u) = 0, \ \forall i \in I_n^\ell, \ n = 1,...,N \label{comp1} \\
&\pi^\top(\sum_{n=1}^N\sum_{i \in I_n^\ell} A_n^i \lambda_{ni} - t) = 0, \ n = 1,...,N \label{comp2} \\
&u_n(\sum_{i \in I_n^\ell} \lambda_{ni} - 1) = 0, \ n = 1,...,N \label{comp3}
\end{align}
Note that $\lambda^{k+1}_{ni}(c_n^i - A_n^{i^\top} \pi^{k+1}_n - u^{k+1}_n) = 0$ from the optimality conditions in (\ref{OPT:ARDMn}). Thus, plugging $\pi^{k+1}, u^{k+1}_n, \lambda^{k+1}_{ni}$ into (\ref{comp1}), we get for each $n$:
\begin{align*}
\abs{\lambda^{k+1}_{ni}(c_n^i - A_n^{i^\top}\pi^{k+1} - u^{k+1}_n)} &\begin{multlined}[t][.6\textwidth]= \Big| \lambda^{k+1}_{ni}(c_n^i - A_n^{i^\top}\pi^{k+1} - u^{k+1}_n) - \lambda^{k+1}_{ni}(c_n^i - A_n^{i^\top}\pi^{k+1}_n - u^{k+1}_n) \Big|\end{multlined} \\
&= \abs{\lambda^{k+1}_{ni}A_n^{i^\top}(\pi^{k+1}_n - \pi^{k+1})} \\
&\leq \abs{\lambda^{k+1}_{ni}} \ \norm{A_n^i} \epsilon_d, \ \forall i \in I_n^\ell 
\end{align*}
Summing over $i \in I_n^\ell$, we get
\begin{align*}
\sum_{i \in I_n^\ell} \abs{\lambda^{k+1}_{ni}(c_n^i - A_n^{i^\top}\pi^{k+1} - u^{k+1}_n)} &\leq \sum_{i \in I_n^\ell} \lambda_{ni}^{k+1} \norm{A_n^i}\epsilon_d \\
&\leq \max_{i \in I_n^\ell} \{\norm{A_n^i}\}\epsilon_d \\
&\leq \norm{A_n}_FL_n\epsilon_d 
\end{align*}
The first inequality follows because $\abs{\lambda_{ni}^{k+1}} = \lambda_{ni}^{k+1}$ and the second inequality holds because $\sum_{i \in I_n^\ell} \lambda_{ni}^{k+1} = 1$. Summing over $n$ gives us
\begin{align}
&\abs{\sum_{n=1}^N\sum_{i \in I_n^\ell}c_n^i \lambda^{k+1}_{ni} - \sum_{n=1}^N\sum_{i \in I_n^\ell}\lambda^{k+1}_{ni}(A_n^{i^\top} \pi^{k+1} + u^{k+1}_n)} \leq \sum_{n=1}^N \norm{A_n}_FL_n\epsilon_d \nonumber \\
\Rightarrow &\abs{\hat{z}_{RMP} - \bigg[\sum_{n=1}^N\sum_{i \in I_n^\ell}\lambda^{k+1}_{ni}(A_n^{i^\top} \pi^{k+1} + u^{k+1}_n)\bigg]} \leq \sum_{n=1}^N \norm{A_n}_FL_n\epsilon_d \label{primalbound}
\end{align}

Moreover, using Theorem \ref{thm1}, (\ref{comp2}) becomes
\begin{align*}
\abs{\pi^{{k+1}^\top}(\sum_{i=1}^N\sum_{i \in I_n^\ell} A_n^i \lambda_{ni}^{k+1} - t)} \leq \norm{\pi^{k+1}} N\epsilon_p \leq mM N\epsilon_p
\end{align*}

Since $u^{k+1}_n(\sum_{i \in I_n^\ell} \lambda^{k+1}_{ni} - 1) = 0$ is satisfied for all $n$, we have
\begin{align}\label{dualbound}
&\abs{\pi^{{k+1}^\top}(\sum_{n=1}^N\sum_{i \in I_n^\ell}A_n^i \lambda^{k+1}_{ni}  - t) + \sum_{n=1}^N u^{k+1}_n(\sum_{i \in I_n^\ell} \lambda^{k+1}_{ni} - 1)} \leq mM N\epsilon_p \nonumber \\
\Rightarrow &\abs{\bigg[\sum_{n=1}^N\sum_{i \in I_n^\ell}\lambda^{k+1}_{ni}(A_n^{i^\top} \pi^{k+1} + u^{k+1}_n)\bigg] - \bigg[\pi^{{k+1}^\top}t+ \sum_{n=1}^Nu^{k+1}_n \bigg]} \leq mM N\epsilon_p \nonumber \\
\Rightarrow &\abs{\bigg[\sum_{n=1}^N\sum_{i \in I_n^\ell}\lambda^{k+1}_{ni}(A_n^{i^\top} \pi^{k+1} + u^{k+1}_n)\bigg] - \hat{z}_{RDM}} \leq mM N\epsilon_p
\end{align}

Adding (\ref{primalbound}) and (\ref{dualbound}) gives us
\begin{align*}
&\abs{\hat{z}_{RMP} - \hat{z}_{RDM}} \leq \epsilon_d\sum_{n=1}^N \norm{A_n}_FL_n  + mM N\epsilon_p
\end{align*}
	
\end{proof}

\begin{theorem}[Optimality Gap]\label{thm2}
CDWD terminates with a solution $\hat{x}$ such that:
\begin{align*}
-\epsilon - \gamma \epsilon_d - mMN\epsilon_p \leq \hat{z}_{RMP} - z_{MP}^*\leq 2\gamma\epsilon_d + mMN\epsilon_p 
\end{align*}
where $\gamma = \sum_{n=1}^N \ \norm{A_n}_FL_n$. 
\end{theorem}
\begin{proof}

By Proposition \ref{prop:gap_dual},
\begin{align}\label{lemma_bound}
-\epsilon \leq \hat{z}_{RDM} - z_{DM}^*\leq \epsilon_d\sum_{n=1}^N \norm{A_n}_FL_n
\end{align}
and by Proposition \ref{prop:restricted_duality_gap}: 
\begin{align}\label{eq:restricted_duality_gap}
-\epsilon_d\sum_{n=1}^N \norm{A_n}_FL_n  - mM N\epsilon_p \leq \hat{z}_{RMP} - \hat{z}_{RDM} \leq \epsilon_d\sum_{n=1}^N \norm{A_n}_FL_n + mM N\epsilon_p
\end{align}

Letting $\gamma = \sum_{n=1}^N \ \norm{A_n}_FL_n$ and adding (\ref{lemma_bound}) and (\ref{eq:restricted_duality_gap}), we get
\begin{align*}
&-\epsilon - \gamma \epsilon_d - mMN\epsilon_p \leq \hat{z}_{RMP} - z_{DM}^* \leq 2\gamma\epsilon_d + mMN\epsilon_p \\
\Rightarrow &-\epsilon - \gamma \epsilon_d - mMN\epsilon_p \leq \hat{z}_{RMP} - z_{MP}^*\leq 2\gamma\epsilon_d + mMN\epsilon_p 
\end{align*}
where the second line of inequalities follows from strong duality, i.e $z_{DM}^* = z_{MP}^*$.
\end{proof}

\section{Computational Experiments}\label{sec:comp_results}

In this section, we present preliminary computational experiments where we solve cutting stock instances from the literature and synthetic instances. CDWD and DWD are implemented in Python using a Message Passing Interface package called mpi4py \cite{dalcin2011parallel,dalcin2008mpi,dalcin2005mpi} and Gurobi is used as the backbone solver. In DWD, only the subproblems are parallelized. The quadratic programs resulting from the augmented Lagrangian \eqref{OPT:ARDMn} are solved using the barrier method. The master problem in DWD and all subproblems are solved using Gurobi's concurrent method. The cutting stock instances were run on a 3.6 GHz Linux machine with 6 cores and 2 threads per core, and the synthetic instances were run on a 3.00 GHz Amazon server running on Linux, with 16 cores and 2 threads per core. We however limit ourselves to only using one thread per core to avoid the overhead of hyperthreading.

\subsection{ADMM Parameters}
To pick the penalty parameter $\rho$, we follow the guidelines provided in \cite{boyd}, where we dynamically adjust $\rho$ according to the primal and dual residuals, so that they are a factor of $\mu$ away from each other, by either multiplying or dividing $\rho$ by positive scalars $\tau^{inc}$ or $\tau^{dec}$, respectively. At the end of iteration $k$, we update $\rho$ as follows:
\begin{align*}
\rho^{k+1} \leftarrow \begin{cases} 
\tau^{inc}\rho^k \ &\text{if} \ \norm{r_d} > \mu \norm{r_p} \\ \frac{\rho^k}{\tau^{dec}} \ &\text{if} \ \norm{r_p} > \mu \norm{r_d} \\ \rho^k \ &\text{otherwise}
\end{cases}
\end{align*}
Intuitively, increasing $\rho$ would put more weight on the terms $\norm{\pi - \pi_n}^2$, thus reducing dual feasibility violation, and decreasing $\rho$ would put more weight on dual optimality, reducing primal feasibility violation. 

We also dynamically adjust the tolerances, similar to \cite{gondzio3}. We solve the pricing subproblems with increasingly accurate dual solutions, where we first solve CDWD with high ADMM tolerances, then divide the tolerances by 10. We repeat the process until we reach desired target tolerances. This significantly reduces computation time, as we reduce the number of times we solve ADMM to high accuracy. We note that the threshold used to add a new column depends on the dual tolerance $\epsilon_d$ (see Section \ref{subsec:cdwd_alg} and Lemma \ref{lem:sep_bounds1}). When starting with high dual tolerances at the beginning of the algorithm, the bound derived in Lemma \ref{lem:sep_bounds1} can be too loose, making the condition to add new columns too harsh. We found it computationally beneficial to be more lenient in adding columns by setting the threshold according to the target dual tolerance and not adapting the threshold according to the current dual tolerance used in solving CDWD.

In our experiments, we notice achieving desired dual tolerances to be harder than reducing feasibility violation. To this end, we pick $\tau^{dec} < \tau^{inc}$ to drive $\rho$ up more easily (thus reducing dual feasibility violation). Moreover, feasibility violations tend to be low, so we pick slightly higher primal feasibility tolerances.

As in classical Dantzig-Wolfe Decomposition (DWD), the RMP might be infeasible if the starting set of extreme points is too small. We circumvent this by adding upper and lower bounds $-M_n \leq \pi_n \leq M_n$ for each block $n$ to ensure a bounded dual problem. This is equivalent to solving the RMP using a big-M method as discussed in Lemma \ref{lem:bound_dual_components}. Computing $M_n$ exactly, however, is prohibitive, and we simply set $M_n = 10\norm{c_n}$ for each block $n$.

We report the optimality gap, computed as $\frac{|\hat{z}_{RMP} - z^*_{MP}|}{|z^*_{MP}|}$, where $\hat{z}_{RMP}$ is the objective value of the RMP evaluated at the recovered primal solution of CDWD and $z^*_{MP}$ is the optimal objective value of the instance. We note that the cutting stock model can only be solved using a column generation algorithm, and many of the larger synthetic instances could not be solved by Gurobi as the process is killed due to the size of the problem. We thus use the objective value obtained from DWD for $z^*_{MP}$.  To compute the feasibility violation, we compute the differences between the left and right-hand sides for each violated linking constraint, normalize by the right-hand side, and report the maximum value. We provide more details as we present the instances.

\subsection{Cutting Stock Problem}
 In the cutting stock problem (CSP) with multiple stock lengths, we are given $K$ types of rolls of different lengths $L_k$ and costs $c_k$, and $P$ pieces with demands $d_p$ and lengths $\ell_p$. Indexing the rolls of each type by $\{1,...,N_k\}$, where $N_k$ is an upperbound on the number of stocks of type $k$ needed, the CSP with multiple stock lengths can be modeled as:
\begin{mini!}
	{x, y}
	{\sum_{k=1}^K \sum_{n=1}^{N_k} c_k y_{kn} \nonumber}
	{\label{OPT:CSP_original}}
	{}
	\addConstraint{\sum_{k=1}^K\sum_{n=1}^{N_k} x_{knp}}{\geq d_p, \label{OPT:CSP.link}}{\ p = 1,\dots,P}
	\addConstraint{\sum_{p=1}^P \ell_p x_{knp}}{\leq L_k, \label{OPT:CSP.block1}}{\ n = 1, \dots, N_k, \ k = 1,...,K}
	\addConstraint{y_{kn}}{\in \{0, 1\},}{\ \forall n, \ \forall k}
	\addConstraint{x_{knp}}{\in \mathbb{Z}_+, \label{OPT:CSP.block2}}{\ \forall p, \ \forall n,\ \forall k}
\end{mini!}
where $y_{kn}$ is one if roll $n \in \{1,\dots,N_k\}$ of type $k$ is used, and $x_{knp}$ is the number of pieces of type $p$ cut from roll $n \in \{1,\dots,N_k\}$ of type $k$. Constraints \eqref{OPT:CSP.link} correspond to the linking constraints, ensuring demand satisfaction. For each $n \in \{1,\dots,N_k\}$ and $k \in \{ 1,\dots,K\}$, constraints \eqref{OPT:CSP.block1}-\eqref{OPT:CSP.block2} correspond to a block's constraints, ensuring the sum of the lengths of the pieces cut from a roll do not exceed the roll's length. There are then $\sum_{k=1}^K N_k$ blocks.

\subsubsection{Reformulation}
Given extreme points $\{x^i_{kn}, y^i_{kn}\}_{i \in I_{kn}}$ for each block $(k, n)$, the master problem can be written as
\begin{mini*}
	{x, y}
	{\sum_{k=1}^K \sum_{n=1}^{N_k} \sum_{i \in I_{kn}} c_k \lambda_{kni}y_{kn}^i \nonumber}
	{}
	{}
	\addConstraint{\sum_{k=1}^K\sum_{n=1}^{N_k} \sum_{i \in I_{kn}} \lambda_{kni} x_{knp}^i}{\geq d_p,}{\ p = 1,\dots,P}
	\addConstraint{\sum_{i \in I_{kn}} \lambda_{kni}}{= 1,}{\ n = 1, \dots, N_k, \ k = 1,...,K}
	\addConstraint{\lambda_{kni}}{\geq 0}{\ \forall i \in I_{kn}, \ \forall n, \ \forall k.}
\end{mini*}
which is a relaxation of the original CSP. For a fixed stock type $k$, the resulting subproblems for all $n \in \{1,...,N_k\}$ are the same and thus return the same column. As in \cite{gondzio1}, we aggregate variables such that $\lambda_{ki} = \sum_{n=1}^{N_k} \lambda_{kni}$ for all $i \in I_k$, so that we have one subproblem per stock type $k$. The resulting master problem is
\begin{mini}
	{x, y}
	{\sum_{k=1}^K \sum_{i \in I_k} c_k \lambda_{ki}y_{k}^i \nonumber}
	{\label{OPT:CSP} \tag{$CSP$}}
	{}
	\addConstraint{\sum_{k=1}^K \sum_{i \in I_k} \lambda_{ki} x_{kp}^i}{\geq d_p,}{\ p = 1,\dots,P}
	\addConstraint{\sum_{i \in I_k} \lambda_{ki}}{= N_k,}{\ k = 1,...,K}
	\addConstraint{\lambda_{ki}}{\geq 0}{\ \forall i \in I_k, \ \forall k,}
\end{mini}
and the resulting subproblem for block $k$ is
\begin{mini}
	{x}
	{c_k - \sum_{p=1}^P \pi_p^* x_{p}}
	{\label{OPT:CSP.subproblem} \tag{$CSP_k$}}
	{}
	\addConstraint{\sum_{p=1}^P \ell_p x_p}{\leq L_k}
	\addConstraint{x_p}{\in \mathbb{Z}_+,}{\ p = 1,\dots,P,}
\end{mini}
where $\pi^*$ corresponds to the optimal dual variables associated with the linking constraints. To add a new column, the solution $y_k^i = 1$ and $x^i_{kp} = x^*_p$ is added as a new extreme point $i$, where $x^*$ is the optimal solution of \eqref{OPT:CSP.subproblem}. The extreme points $x^i_k$ are feasible cutting patterns of a stock of type $k$ and variables $\lambda$ are selecting patterns such that demands are satisfied in \eqref{OPT:CSP}. 

Finally, we note that since the solution $x_p = 0$ for all $p$ is feasible in \eqref{OPT:CSP.subproblem} and has cost zero, we can relax the convexity constraints to $\sum_{i \in I_k} \lambda_{ki} \leq N_k$ and since $N_k$ is an upperbound on the number of stocks of type $k$ needed, we can omit the convexity constraints in \eqref{OPT:CSP} (see \cite{gondzio1} for more details). This omission is accounted for in \eqref{OPT:CSP.subproblem}.

\subsubsection{Computational Results}
We solve CSP instances with multiple stock lengths obtained from the CaPaD library \cite{belov2002cutting}(\url{http://www.math.tu-dresden.de/~capad/}). We solve instances with 4 and 5 stock lengths, and approximately 40, 150, 200, 300, and 400 items. For the ADMM parameters, we pick $\mu = 50$, $\tau^{inc} = 2$, $\tau^{dec} = 1.5$ and $\rho^0 = 100$. For tolerances, we start at $\epsilon_p = 5$ and $\epsilon_d = 50$, and end at target tolerances $\epsilon_p = 5 \times 10^{-2}$ and $\epsilon_d = 5 \times 10^{-3}$.

In Table \ref{tab:csp}, we include the range of the number of items $P$ across the instances, the number of stock lengths $K$, the geometric means (GM) of the optimality gaps, the feasibility violations and the runtimes, and the medians of the runtimes. Given a recovered solution $\hat x$, the feasibility violation is computed as
\begin{align*}
	\max_{p}\left\{ \frac{d_p - \sum_{k=1}^K \hat x_{kp}}{d_p}, 0 \right \}.
\end{align*}

\begin{table}
	\centering
	\caption{CSP Results} \label{tab:csp}
	{\footnotesize
	\begin{tabular}{C{1.1cm}C{0.5cm}C{2cm}C{2cm}C{2cm}C{2cm}C{2cm}C{2cm}}
		\toprule[1pt]\midrule[0.3pt]
		$P$ Range &$K$ &Optimality Gap &Feasibility Gap  & CDWD Time - GM (sec)  &DWD Time - GM (sec)  &CDWD Time - Median (sec)  &DWD Time - Median (sec) \\ 
		\hline
		$\interval{38}{40}$ 		&5&1.00e-02	&2.46e-04	&1.3	&2.19		&1.15		&2.23\\
		$\interval{146}{150}$	&4	&1.00e-02	&8.85e-03	&12.63	&17.53		&11.39		&19.06\\
		$\interval{194}{200}$ 	&4	&1.00e-02	&8.90e-03	&22.32	&29.99		&21.21		&31.36\\
		$\interval{289}{299}$	&4	&9.60e-03	&9.43e-03	&64.52	&58.23		&58.66		&60.84\\
		$\interval{385}{395}$ 	&4	&1.03e-02	&9.38e-03	&138.83	&101.5		&116.03		&103.26 \\
		\midrule[0.3pt]\bottomrule[1pt]
	\end{tabular}
}
\end{table}

CDWD recovers high quality solutions with low optimality gaps and low feasibility violations. The geometric mean of optimality gaps is about $10^{-2}$ and the geometric mean of feasibility violations is close to $10^{-4}$ for all five sets of experiments. We note that CDWD runs faster than DWD in most instances in the first three sets of experiments, with the geometric means of CDWD's runtimes being about $40\%$, $28\%$ and $26\%$ lower than DWD's geometric means, and the medians of CDWD's runtimes being about $48\%$, $40\%$, $32\%$ lower than DWD's medians, respectively. For instances with about 200 items, the geometric means are close, with DWD's being slightly lower, but CDWD runs faster in over $50\%$ of the instances. With about 400 items, DWD runs faster in most instances, with a geometric mean that is about $27\%$ lower and a median that is about $11\%$ lower. As we increase the number of items, the difference between the runtimes decreases, where DWD eventually runs faster in most instances with $P$ close to 400. We note that the number of linking constraints is equal to the number of items, thus increasing linearly as the CSP instances get larger. As the number of linking constraints gets larger, ADMM convergence starts to slow down, explaining our observations. To this end, we turn to synthetic instances to more easily perform sensitivity analysis with respect to the number of blocks, variables and linking constraints.

\subsection{Synthetic Instances}
\subsubsection{Instance Generation}
The synthetic instances are of the form
\begin{align*}
\min& \qquad \sum_{n=1}^N c_n^\top x_n \\
\text{s.t.}& \qquad \sum_{n=1}^NA_nx_n \geq t  \\
& \qquad B_nx_n \leq b_n, \ \forall n = 1,...,N \\
& \qquad 0 \leq x_n \leq u_n, \ \forall n = 1,...,N
\end{align*}
where the coefficients of the matrices $A_n$ and $B_n$ are from the discrete uniform distribution $\mathcal{U}\{-10,20\}$, and the components of the cost vector are from $\mathcal{U}\{-10,30\}$. Let $\ell_i$ be the sum of the entries in row $i$ of the linking constraints, i.e $\ell_i = \sum_{n,j} (A_n)_{ij}$, where $(A_n)_{ij}$ is component $(i,j)$ of $A_n$; similarly let $\beta^n_i$ be the sum of the entries of row $i$ of $B_n$. The vectors $t$ and $b_n$ were generated according to the sum of each row of the constraint matrix. We construct component $i$ of $t$ as follows:	
\begin{align*}
\begin{cases} t_i \sim \mathcal{U}\{2\ell_i, 3\ell_i\}, \ &\text{if} \ \ell_i > 0 \\
						t_i \sim \mathcal{U}\{3\ell_i, 2\ell_i\}, \ &\text{if} \ \ell_i < 0 \\
						t_i = 0, \ &\text{if} \ \ell_i = 0
\end{cases} i = 1,...,m
\end{align*}
where $m$ is the number of linking constraints. Similarly, component $i$ of $b_n$ is constructed as
\begin{align*}
\qquad \begin{cases} (b_n)_i \sim \mathcal{U}\{2\beta^n_i, 3\beta^n_i\} \ &\text{if} \ \beta^n_i > 0 \\
						(b_n)_i \sim \mathcal{U}\{3\beta^n_i, 2\beta^n_i\}, \ &\text{if} \ \beta^n_i < 0 \\
						(b_n)_i = 0, \ &\text{if} \ \beta^n_i = 0
\end{cases} i = 1,...,m_n
\end{align*}
where $m_n$ is the number of constraints in block $n$. Moreover, to ensure a bounded region, we add upper and lower bounds to the variables, where $u_n = 30$ for all $n$.

\subsubsection{Computational Results}
We perform four sets of experiments, each set involving 1, 2, 5, and 10 linking constraints. The results are reported in Tables \ref{tab:1link}-\ref{tab:10link}. For each set of experiment, we vary the number of blocks $N \in \{2, 4, 8, 10, 15\}$ and the total number of variables across all blocks $n_v \in \{100, 1000, 5000, 10000, \\ 20000, 25000, 50000, 100000\}$. We define $m_n$ to be the number of block constraints. For simplicity, each block has approximately the same number of variables. We limit the number of blocks to 15 to avoid the overhead of hyperthreading. We note that we have also solved the instances as a single LP using Gurobi but omit the results, as Gurobi's runtimes, although lower for the smaller instances, quickly become greater than CDWD and DWD as the instances get larger. Moreover, Gurobi failed to solve instances with $n_v \geq 50000$.
The feasibility violations are computed as
\begin{align*}
	\max_i\left\{ \frac{t_i - \sum_{n,j}(A_n)_{ij}(\hat{x}_n)_j}{|t_i|}, 0 \right \}.
\end{align*}
For the ADMM parameters, we pick $\mu = 100$, $\tau^{inc} = 2$, $\tau^{dec} = 1.5$ and $\rho^0 = 100$. For tolerances, we start at $\epsilon_p = 50$ and $\epsilon_d = 50$, and end at target tolerances $\epsilon_p = 5 \times 10^{-2}$ and $\epsilon_d = 5 \times 10^{-3}$.

Across all instances, feasibility violations are very close to 0, with the largest violation being in the order of $10^{-4}$; optimality gaps are also very low, with most instances hovering between the order of $10^{-6}$ and $10^{-2}$, with only eight instances reaching optimality gaps as high as $10^{-2}$.

\begin{figure}[t!]
	\centering
    \begin{subfigure}[c]{0.45\textwidth}
        \includegraphics[scale=0.5]{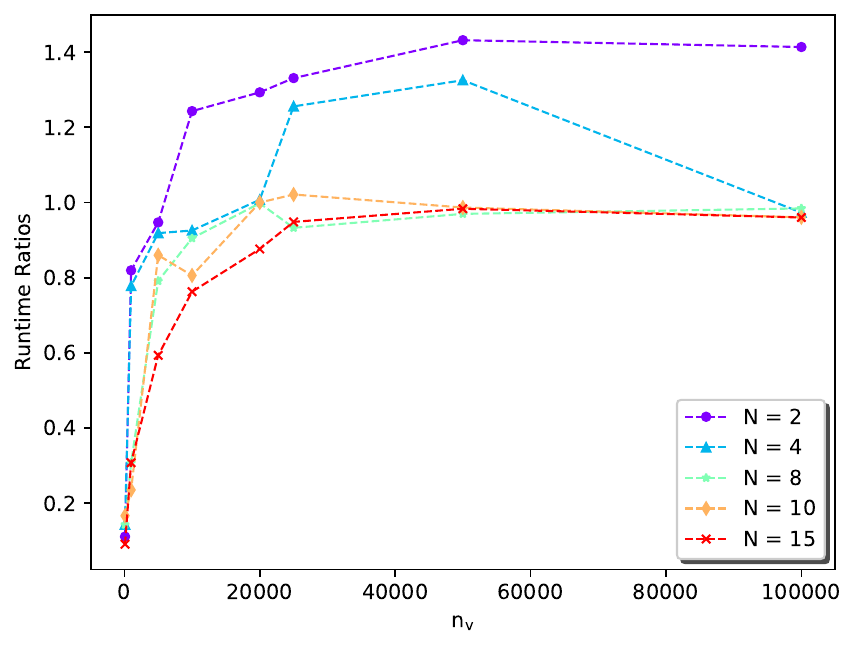}
        \caption{$m=1$}\label{fig:1link}
    \end{subfigure} \hspace{5ex}
    \begin{subfigure}[c]{0.45\textwidth}
        \includegraphics[scale=0.5]{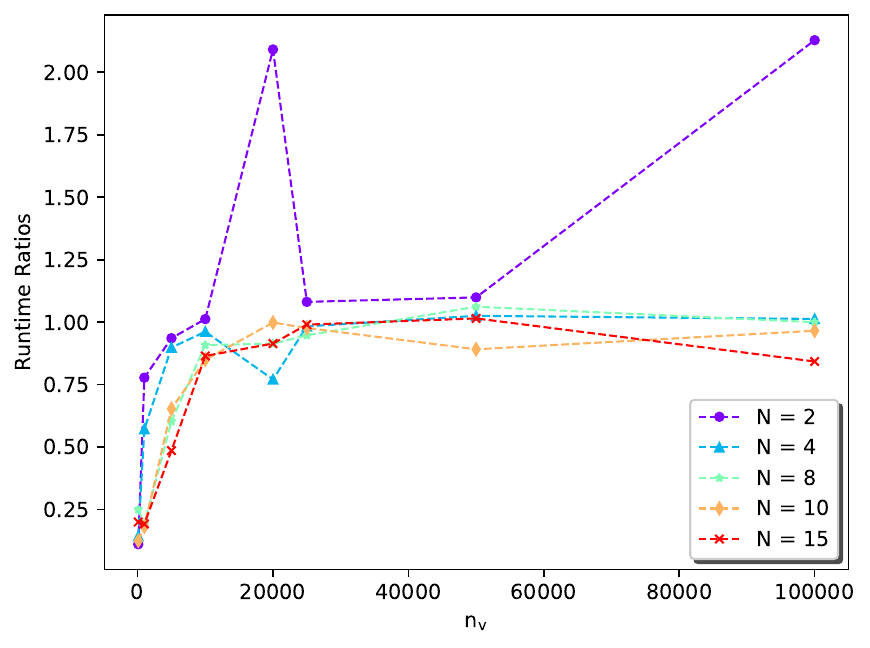}
        \caption{$m=2$}\label{fig:2link}
	\end{subfigure} \\

	\begin{subfigure}[c]{0.45\textwidth}
        \includegraphics[scale=0.5]{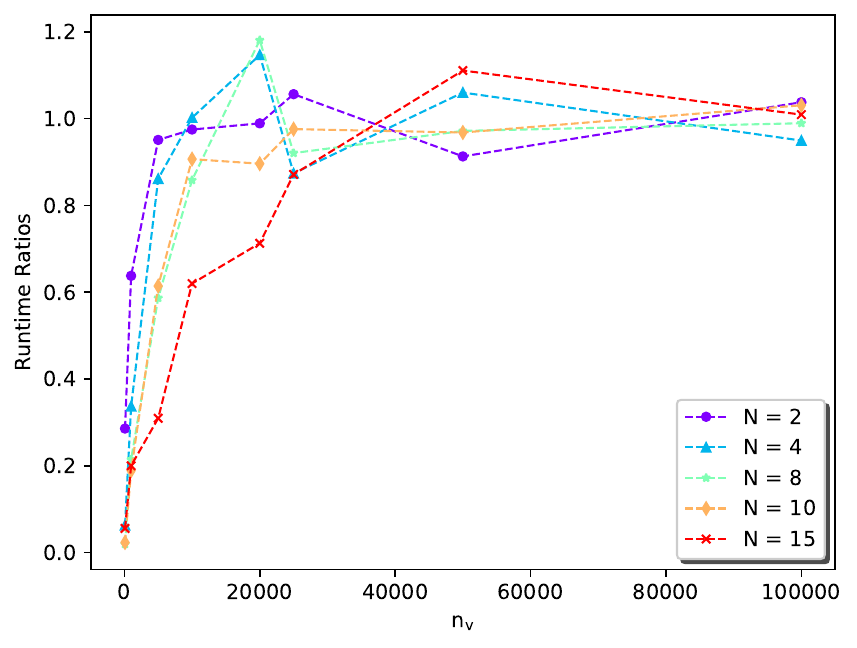}
        \caption{$m=5$}\label{fig:3link}
	\end{subfigure} \hspace{5ex}
	\begin{subfigure}[c]{0.45\textwidth}
        \includegraphics[scale=0.5]{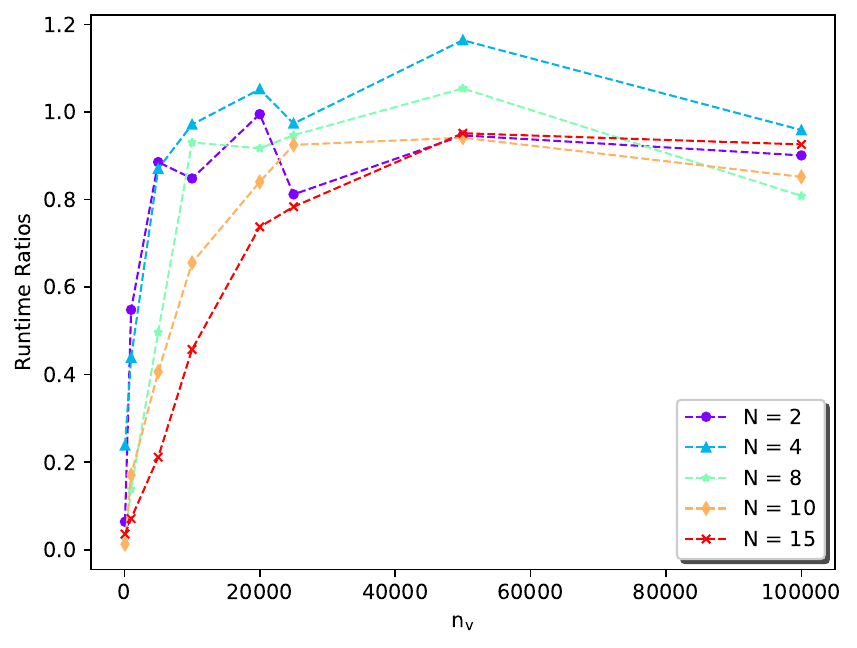}
        \caption{$m=10$}\label{fig:4link}
    \end{subfigure}
	\caption{DWD to CDWD Runtime Ratios}
	\label{fig:ratios} 
\end{figure}

For smaller instances, DWD runs faster than CDWD. For a fixed number of linking constraints and blocks, the gap between the runtimes closes as the number of variables increases, with CDWD eventually running faster in many instances with 25000 variables or more, especially for smaller values of $m$ and $N$.  For example, with five linking constraints and ten blocks, DWD is over an order of magnitude faster for $n_v = 100$ and $n_v = 1000$, and is about 50\% faster for $n_v = 5000$; the runtimes are almost identical for larger values of $n_v$, until CDWD becomes faster at 100000 total variables. We observe this overall trend in many cases, but note that the gap between the runtimes levels out, and further experiments may be required to confirm it. We plot the ratios of DWD to CDWD runtimes in Figure \ref{fig:ratios} (a value under 1 means DWD is faster). The ratio approaches one fast as we increase $n_v$ to 25000, after which CDWD runs faster in more and more instances. Between 50000 and 100000 variables, the ratios either stay close to one, slightly decrease in some cases or slightly increase in others. The ratio also approaches one at a slower rate as we increase $N$ and, to some extent, the number of linking constraints. This is expected as increasing the number of blocks and linking constraints slows down the convergence of ADMM, thus requiring larger instances in terms of number of variables for CDWD to catch up to DWD and potentially run faster.

Although our main focus is dealing with privacy issues and data that cannot be stored in a central location, either for security reasons or physical limitations, we show that the computational price of handling these concerns approaches zero as the number of variables increases, and even see a potential computational benefit in many cases, while sacrificing very little in terms of solution quality.

\begin{table}
\centering
\caption{1 Linking Constraint} \label{tab:1link}
{\footnotesize
\begin{tabular}{C{0.5cm}M{1.1cm}M{1cm}C{2cm}C{2cm}M{1.1cm}M{1.1cm}}
\toprule[1pt]\midrule[0.3pt]
$N$ &$n_v$ &$m_n$ &Optimality Gap &Feasibility Gap  &CDWD Time (sec)  &DWD Time (sec) \\ 
\hline
\multirow{8}{*}{\begin{tabular}[M{2.5cm}]{@{}l@{}}2\end{tabular}}	
&100	&50		&9.89e-05	& 0 	&0.09	&0.01	\\
&1000	&500	&5.07e-06	& 0 	&0.72	&0.59	\\
&5000	&2500	&6.87e-06	& 0 	&32.1	&30.4	\\
&10000	&2500	&2.89e-02	& 0 	&57.6	&71.6	\\
&20000	&2500	&1.07e-02	& 0 	&123.2	&159.3	\\
&25000	&2500	&5.86e-03	& 0 	&155.3	&206.7	\\
&50000	&2500	&9.87e-03	& 0 	&357.6	&512.0	\\
&100000	&2500	&2.41e-02	& 0 	&801.4	&1132.9	\\
\hline
\multirow{8}{*}{\begin{tabular}[M{2.5cm}]{@{}l@{}}4 \end{tabular}}	
&100	&25		&1.71e-04	& 0 	&0.07	&0.01		\\
&1000	&250	&9.97e-05	& 0 	&0.18	&0.14		\\
&5000	&1250	&4.36e-05	&2.35e-06	&6.75	&6.2		\\
&10000	&2500	&1.16e-06	&1.11e-07	&37.6	&34.8		\\
&20000	&2500	&4.60e-06	& 0 	&82.3	&82.9		\\
&25000	&2500	&2.75e-02	& 0 	&82.9	&104.1		\\
&50000	&2500	&2.17e-02	& 0 	&173.7	&230.2	\\
&100000	&2500	&1.05e-05	& 0 	&531.2	&516.9	\\
\hline
\multirow{8}{*}{\begin{tabular}[M{2.5cm}]{@{}l@{}}8 \end{tabular}}	
&100	&12	&2.41e-04	&2.06e-04	&0.07	&0.01		\\
&1000	&125	&2.17e-05	& 0 	&0.13	&0.04		\\
&5000	&625	&5.84e-06	&4.72e-06	&1.82	&1.44		\\
&10000	&1250	&2.89e-05	& 0 	&9.99	&9.04		\\
&20000	&2500	&1.57e-05	&8.13e-09	&50.1	&50.0		\\
&25000	&2500	&7.96e-06	& 0 	&69.8	&65.1		\\
&50000	&2500	&5.76e-06	& 0 	&148.7	&144.2	\\
&100000	&2500	&1.30e-05	&1.92e-07	&316.6	&311.6	\\
\hline
\multirow{8}{*}{\begin{tabular}[M{2.5cm}]{@{}l@{}}10 \end{tabular}}	
&100	&10		&2.81e-05	&1.50e-04	&0.06		&0.01		\\
&1000	&100	&5.37e-05	&1.18e-05	&0.17		&0.04	\\
&5000	&500	&1.48e-05	&2.66e-06	&1.57		&1.35	\\
&10000	&1000	&9.41e-06	&3.32e-06	&9.32		&7.51	\\
&20000	&2000	&4.75e-06	& 0 		&43.4		&43.4		\\
&25000	&2500	&1.48e-05	&7.61e-08	&75.3		&76.9	\\
&50000	&2500	&1.97e-05	& 0 		&164.6		&162.4	\\
&100000	&2500	&1.31e-05	& 0 		&370.6		&356.0	\\
\hline
\multirow{8}{*}{\begin{tabular}[M{2.5cm}]{@{}l@{}}15 \end{tabular}}	
&100	&6	&6.88e-04	&5.39e-04	&0.11	&0.01		\\
&1000	&66	&1.01e-04	& 0 	&0.13	&0.04		\\
&5000	&333	&3.49e-05	& 0 	&0.91	&0.54		\\
&10000	&666	&4.05e-06	& 0 	&3.96	&3.02		\\
&20000	&1333	&9.40e-06	& 0 	&19.4	&17.0		\\
&25000	&1666	&9.55e-06	&5.28e-07	&31.0	&29.4	\\
&50000	&2500	&9.55e-06	& 0 	&120.2	&118.2	\\
&100000	&2500	&5.08e-06	& 0 	&265.8	&255.2	\\
\midrule[0.3pt]\bottomrule[1pt]
\end{tabular}
}
\end{table}
\begin{table}
	\centering
	\caption{2 Linking Constraints} \label{tab:2link}
	{\footnotesize
	\begin{tabular}{C{0.5cm}M{1.1cm}M{1cm}C{2cm}C{2cm}M{1.1cm}M{1.1cm}}
	\toprule[1pt]\midrule[0.3pt]
	$N$ &$n_v$ &$m_n$ &Optimality Gap &Feasibility Gap  &CDWD Time (sec)  &DWD Time (sec) \\ 
	\hline
	\multirow{8}{*}{\begin{tabular}[M{2.5cm}]{@{}l@{}}2\end{tabular}}	
	&100	&50	&3.59e-04	&5.80e-05	&0.09	&0.01		\\
	&1000	&500	&6.61e-05	&9.01e-11	&0.72	&0.56		\\
	&5000	&2500	&8.02e-06	&5.19e-07	&43.5	&40.7		\\
	&10000	&2500	&1.66e-04	& 0 	&89.5	&90.6		\\
	&20000	&2500	&3.37e-02	& 0 	&105.1	&219.8		\\
	&25000	&2500	&5.45e-06	&5.07e-09	&276.1	&298.5		\\
	&50000	&2500	&3.69e-05	&8.52e-09	&690.2	&758.5		\\
	&100000	&2500	&3.03e-02	& 0 	&841.6	&1791.3		\\
	\hline
	\multirow{8}{*}{\begin{tabular}[M{2.5cm}]{@{}l@{}}4 \end{tabular}}	
	&100	&25	&2.40e-03	&1.55e-04	&0.07	&0.01		\\
	&1000	&250	&2.02e-05	&8.76e-06	&0.28	&0.16		\\
	&5000	&1250	&2.48e-05	&1.07e-06	&9.33	&8.38		\\
	&10000	&2500	&2.47e-05	& 0 	&47.0	&45.2		\\
	&20000	&2500	&6.06e-06	& 0 	&102.4	&78.9		\\
	&25000	&2500	&2.66e-04	&1.33e-07	&141.5	&139.1		\\
	&50000	&2500	&4.16e-05	&1.60e-07	&317.0	&325.0		\\
	&100000	&2500	&3.64e-05	&6.39e-08	&793.9	&804.0	 	\\
	\hline
	\multirow{8}{*}{\begin{tabular}[M{2.5cm}]{@{}l@{}}8 \end{tabular}}	
	&100	&12	&2.15e-04	&4.13e-05	&0.04	&0.01		\\
	&1000	&125	&7.60e-05	&5.10e-06	&0.26	&0.05		\\
	&5000	&625	&1.45e-05	&3.02e-06	&2.92	&1.76		\\
	&10000	&1250	&6.34e-06	&9.05e-07	&13.1	&11.9		\\
	&20000	&2500	&2.30e-05	&1.97e-07	&74.1	&67.7		\\
	&25000	&2500	&6.61e-06	&3.54e-07	&100.2	&94.9		\\
	&50000	&2500	&2.32e-05	& 0 	&205.6	&218.3	 	\\
	&100000	&2500	&1.22e-05	&6.43e-08	&447.8	&447.7		\\
	\hline
	\multirow{8}{*}{\begin{tabular}[M{2.5cm}]{@{}l@{}}10 \end{tabular}}	
	&100	&10	&2.18e-03	&5.64e-04	&0.08	&0.01		\\
	&1000	&100	&6.58e-05	& 0 	&0.33	&0.06		\\
	&5000	&500	&1.05e-05	& 0 	&2.39	&1.56		\\
	&10000	&1000	&1.37e-05	&1.44e-06	&9.85	&8.37		\\
	&20000	&2000	&9.50e-06	& 0 	&57.6	&57.5		\\
	&25000	&2500	&2.73e-05	&6.04e-07	&102.1	&99.7		\\
	&50000	&2500	&4.82e-06	& 0 	&251.2	&223.9		\\
	&100000	&2500	&5.69e-06	& 0 	&528.5	&510.2		\\
	\hline
	\multirow{8}{*}{\begin{tabular}[M{2.5cm}]{@{}l@{}}15 \end{tabular}}	
	&100	&6		&7.36e-04	&6.46e-04	&0.1	&0.02	\\
	&1000	&66		&8.07e-05	& 0 		&0.26	&0.05		\\
	&5000	&333	&3.30e-05	&4.21e-06	&1.4	&0.68		\\
	&10000	&666	&1.38e-05	& 0 		&4.36	&3.77		\\
	&20000	&1333	&1.47e-05	&1.16e-06	&24.5	&22.4		\\
	&25000	&1666	&8.19e-06	&1.05e-06	&40.4	&40.0		\\
	&50000	&2500	&4.75e-06	& 0 		&141.6	&143.7		\\
	&100000	&2500	&3.31e-06	&4.27e-08	&411.9	&347.0	\\
	\midrule[0.3pt]\bottomrule[1pt]
	\end{tabular}
	}
\end{table}
\begin{table}
	\centering
	\caption{5 Linking Constraints} \label{tab:5link}
	{\footnotesize
	\begin{tabular}{C{0.5cm}M{1.1cm}M{1cm}C{2cm}C{2cm}M{1.1cm}M{1.1cm}}
	\toprule[1pt]\midrule[0.3pt]
	$N$ &$n_v$ &$m_n$ &Optimality Gap &Feasibility Gap  &CDWD Time (sec)  &DWD Time (sec)\\ 
	\hline
	\multirow{8}{*}{\begin{tabular}[M{2.5cm}]{@{}l@{}}2\end{tabular}}	
	&100	&50		&2.24e-04	& 0 		&0.07	&0.02		\\
	&1000	&500	&1.88e-05	&3.70e-07	&2.07	&1.32		\\
	&5000	&2500	&1.10e-05	&6.97e-07	&113.9	&108.3		\\
	&10000	&2500	&1.69e-05	&1.36e-08	&238.1	&232.1		\\
	&20000	&2500	&3.46e-05	&2.58e-07	&572.6	&566.3		\\
	&25000	&2500	&1.97e-05	&7.51e-08	&766.0	&809.1		\\
	&50000	&2500	&5.35e-05	&6.67e-08	&2040.1	&1862.4	 	\\
	&100000	&2500	&3.45e-05	&7.14e-09	&3202.5	&3323.3	 	\\
	\hline
	\multirow{8}{*}{\begin{tabular}[M{2.5cm}]{@{}l@{}}4 \end{tabular}}	
	&100	&25	&1.47e-04	&7.48e-05	&0.16	&0.01		\\
	&1000	&250	&5.25e-05	& 0 	&0.95	&0.32		\\
	&5000	&1250	&4.53e-05	&2.63e-06	&17.3	&14.9		\\
	&10000	&2500	&3.23e-05	&6.96e-07	&100.4	&100.6		\\
	&20000	&2500	&4.65e-06	&4.45e-07	&213.6	&244.9		\\
	&25000	&2500	&3.08e-05	&2.61e-07	&331.3	&289.5		\\
	&50000	&2500	&2.59e-05	&4.69e-08	&876.7	&929.3		\\
	&100000	&2500	&1.52e-05	&9.98e-11	&2564.8	&2434.8	 	\\
	\hline
	\multirow{8}{*}{\begin{tabular}[M{2.5cm}]{@{}l@{}}8 \end{tabular}}	
	&100	&12	&6.31e-04	&2.11e-06	&0.52	&0.01		\\
	&1000	&125	&1.66e-05	&1.94e-05	&0.42	&0.09		\\
	&5000	&625	&1.04e-05	&2.71e-06	&5.66	&3.32		\\
	&10000	&1250	&5.31e-05	&2.46e-06	&28.2	&24.2		\\
	&20000	&2500	&1.40e-05	&8.65e-08	&127.9	&151.0		\\
	&25000	&2500	&3.43e-05	&6.57e-08	&193.8	&178.4		\\
	&50000	&2500	&9.62e-06	&1.55e-08	&416.0	&404.2		\\
	&100000	&2500	&1.20e-05	&7.08e-08	&1029.3	&1018.3		\\
	\hline
	\multirow{8}{*}{\begin{tabular}[M{2.5cm}]{@{}l@{}}10 \end{tabular}}	
	&100	&10	&4.87e-04	&2.10e-04	&0.43	&0.01	\\
	&1000	&100	&5.54e-06	&2.23e-05	&0.53	&0.1		\\
	&5000	&500	&4.03e-05	&3.94e-06	&4.51	&2.77		\\
	&10000	&1000	&4.91e-06	&1.05e-07	&21.4	&19.4		\\
	&20000	&2000	&1.38e-05	&1.49e-06	&121.2	&108.6		\\
	&25000	&2500	&1.59e-05	&4.81e-07	&222.7	&217.3		\\
	&50000	&2500	&7.98e-06	&1.14e-07	&430.4	&416.6		\\
	&100000	&2500	&4.02e-05	& 0 	&948.9	&978.0		\\
	\hline
	\multirow{8}{*}{\begin{tabular}[M{2.5cm}]{@{}l@{}}15 \end{tabular}}	
	&100	&6	&2.43e-02	& 0 	&0.36	&0.02		\\
	&1000	&66	&2.04e-04	&3.41e-06	&0.3	&0.06		\\
	&5000	&333	&1.14e-05	& 0 	&4.2	&1.3		\\
	&10000	&666	&2.08e-05	&1.34e-06	&11.1	&6.88		\\
	&20000	&1333	&1.00e-05	&1.53e-06	&42.8	&30.5		\\
	&25000	&1666	&9.55e-06	&5.05e-07	&85.5	&74.5		\\
	&50000	&2500	&7.55e-06	&2.49e-07	&300.6	&333.9		\\
	&100000	&2500	&7.63e-06	&7.41e-08	&813.5	&820.9		\\
	\midrule[0.3pt]\bottomrule[1pt]
	\end{tabular}
	}
\end{table}
\begin{table}
	\centering
	\caption{10 Linking Constraints} \label{tab:10link}
	{\footnotesize
	\begin{tabular}{C{0.5cm}M{1.1cm}M{1cm}C{2cm}C{2cm}M{1.1cm}M{1.1cm}}
	\toprule[1pt]\midrule[0.3pt]
	$N$ &$n_v$ &$m_n$ &Optimality Gap &Feasibility Gap  &CDWD Time (sec)  &DWD Time (sec)\\ 
	\hline
	\multirow{8}{*}{\begin{tabular}[M{2.5cm}]{@{}l@{}}2\end{tabular}}	
	&100	&50	&4.87e-05	&4.82e-05	&0.47	&0.03	\\
	&1000	&500	&6.07e-05	&2.00e-06	&6.04	&3.31	\\
	&5000	&2500	&2.44e-05	&6.34e-07	&255.0	&225.8	\\
	&10000	&2500	&4.23e-05	&3.54e-07	&568.5	&482.1	\\
	&20000	&2500	&2.11e-05	&6.73e-08	&1339.1	&1332.4		\\
	&25000	&2500	&2.05e-05	&2.14e-07	&2334.7	&1895.1		\\
	&50000	&2500	&2.74e-05	&8.39e-09	&5412.5	&5122.9		\\
	&100000	&2500	&4.91e-05	&5.20e-09	&15681.6	&14124.8	\\
	\hline
	\multirow{8}{*}{\begin{tabular}[M{2.5cm}]{@{}l@{}}4 \end{tabular}}	
	&100	&25	&3.50e-03	&6.49e-05	&0.21	&0.05		\\
	&1000	&250	&1.93e-04	&8.00e-06	&1.71	&0.75		\\
	&5000	&1250	&3.53e-05	&5.87e-07	&53.9	&46.9		\\
	&10000	&2500	&3.45e-05	&8.81e-07	&359.1	&348.8		\\
	&20000	&2500	&1.60e-05	&2.18e-07	&631.4	&664.2		\\
	&25000	&2500	&9.95e-06	&8.11e-08	&879.6	&856.4		\\
	&50000	&2500	&1.98e-05	&1.50e-07	&1820.5	&2118.9	 	\\
	&100000	&2500	&1.61e-05	&5.15e-08	&5044.7	&4836.2		\\
	\hline
	\multirow{8}{*}{\begin{tabular}[M{2.5cm}]{@{}l@{}}8 \end{tabular}}	
	&100	&12	&2.64e-04	&2.08e-04	&0.43	&0.01		\\
	&1000	&125	&3.47e-05	&1.55e-05	&1.23	&0.17		\\
	&5000	&625	&1.50e-05	&1.93e-06	&11.9	&5.91		\\
	&10000	&1250	&1.38e-05	&1.16e-06	&50.2	&46.7		\\
	&20000	&2500	&1.67e-05	&4.49e-07	&329.8	&302.4		\\
	&25000	&2500	&3.14e-05	&2.60e-07	&372.9	&353.1		\\
	&50000	&2500	&1.20e-05	&3.54e-08	&867.2	&913.8	 	\\
	&100000	&2500	&9.61e-06	&1.11e-07	&2659.3	&2149.4		\\
	\hline
	\multirow{8}{*}{\begin{tabular}[M{2.5cm}]{@{}l@{}}10 \end{tabular}}	
	&100	&10	&5.71e-05	&1.11e-04	&0.8	&0.01		\\
	&1000	&100	&4.86e-05	&1.79e-05	&0.82	&0.14		\\
	&5000	&500	&6.93e-05	&3.68e-06	&11.1	&4.51	\\
	&10000	&1000	&1.61e-05	&1.49e-06	&44.1	&28.9		\\
	&20000	&2000	&2.27e-05	&6.12e-07	&209.1	&175.7		\\
	&25000	&2500	&1.19e-05	&3.32e-07	&412.5	&381.5		\\
	&50000	&2500	&6.69e-06	&2.76e-07	&845.3	&795.6	 	\\
	&100000	&2500	&7.34e-06	&4.53e-08	&2601.6	&2215.6		\\
	\hline
	\multirow{8}{*}{\begin{tabular}[M{2.5cm}]{@{}l@{}}15 \end{tabular}}	
	&100	&6	&1.22e-05	&5.10e-04	&0.56	&0.02		\\
	&1000	&66	&8.78e-05	&1.85e-05	&1.55	&0.11		\\
	&5000	&333	&2.43e-05	&2.58e-06	&8.14	&1.72		\\
	&10000	&666	&1.78e-05	&9.22e-07	&22.3	&10.2		\\
	&20000	&1333	&2.12e-05	&1.27e-06	&99.4	&73.3		\\
	&25000	&1666	&1.63e-05	&3.63e-07	&163.9	&128.4		\\
	&50000	&2500	&9.94e-06	&7.10e-07	&559.7	&532.5		\\
	&100000	&2500	&7.16e-06	&5.97e-08	&1176.1	&1088.9		\\
	\midrule[0.3pt]\bottomrule[1pt]
	\end{tabular}
	}
\end{table}

\subsection{Parallel Efficiency and Scalability}
To measure how our Python implementation scales as we increase the number of blocks and available cores, we use two common metrics as in \cite{munguia2018alternating}. The first one measures the speedup gained by using the available cores. The second metric measures core utilization and time lost in communication and synchronization. We compute the two metrics for instances with 9000, 18000 and 36000 total variables. For each set of instances, we experiment with 5, 10, 20, 36 and 72 blocks. As before, each block contains the same number of variables. Note that in these experiments, we potentially pay the price of hyperthreading. The main objective is to showcase the slowdowns that can be caused by idle cores, indicating potential benefits in an asynchronous implementation of our algorithm.\\

\begin{figure}
\caption{Ratio of runtimes between serial and parallel implementations}
\label{fig:speedups}
\begin{center}
\includegraphics[scale=0.25]{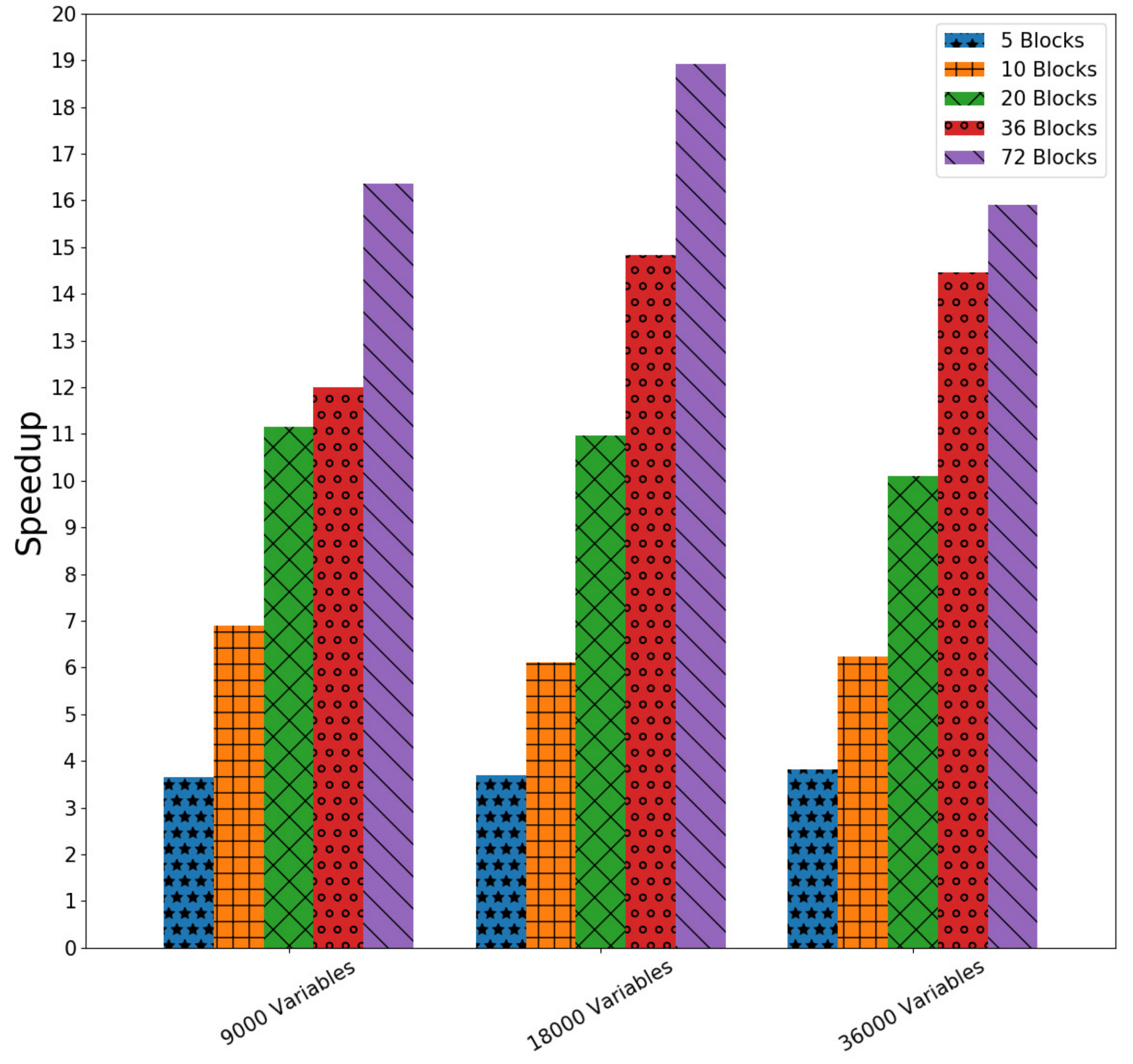}
\end{center} 
\end{figure}

\noindent \textbf{Parallel Speedup} Let $t_p$ be the time it takes for CDWD to terminate using $p$ cores. We compute the ratio $\frac{t_1}{t_p}$ for each experiment and report results in Figure \ref{fig:speedups}. We observe similar trends for different number of total variables. The computational gain from parallelizing decreases as we increase the number of blocks. This is mainly due to cores sitting idle, waiting on other processes to finish, as well as communication overhead increasing with the number of cores used. This is confirmed by our analysis on core utilization.
\\

\begin{figure}
\caption{Average core utilization}
\label{fig:utilization}
\begin{center}
\includegraphics[scale=0.25]{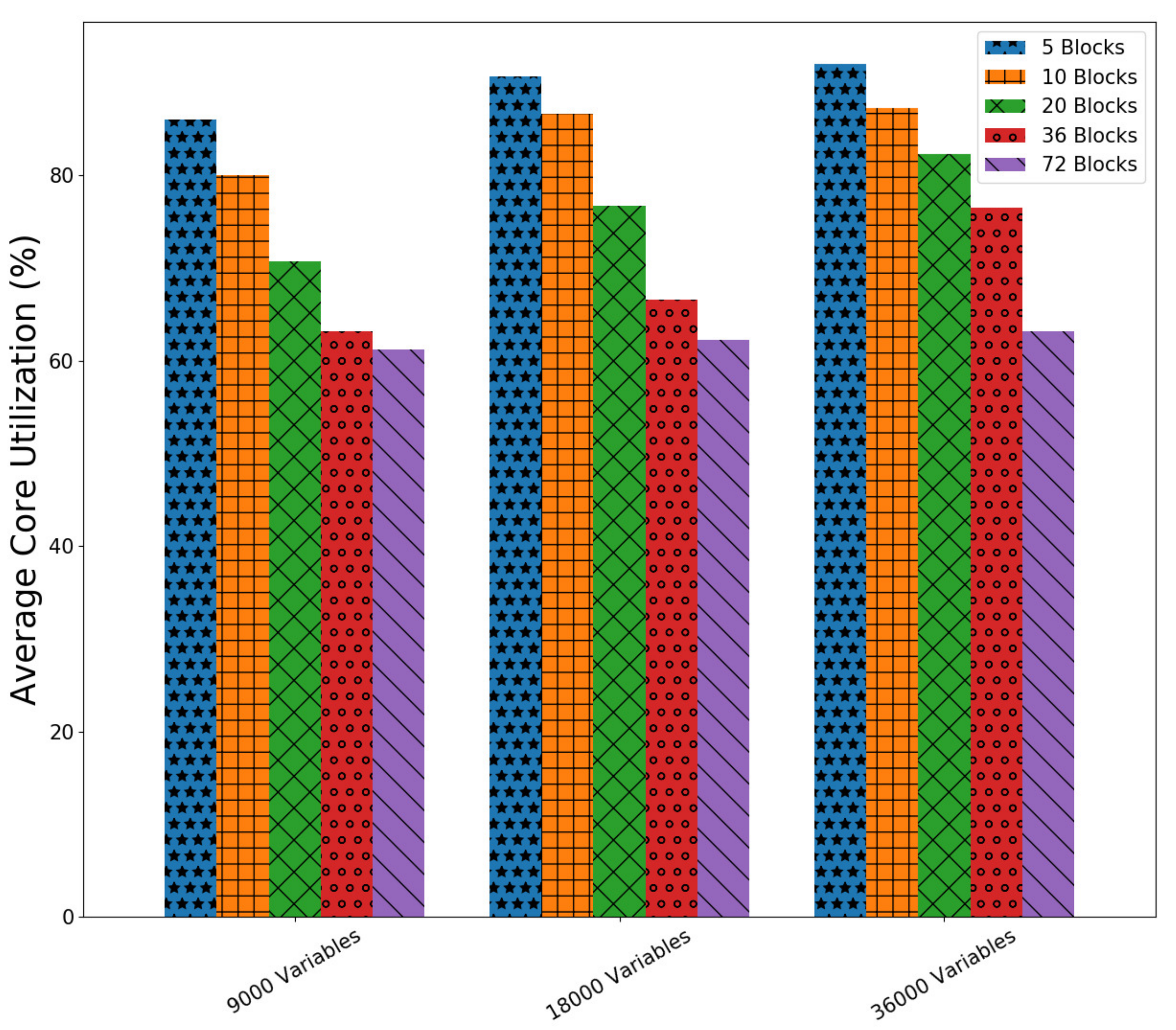}
\end{center} 
\end{figure}

\noindent \textbf{Core Utilization} To estimate core utilization, we measure total time spent doing useful computations, communication time, and synchronization time where a core is sitting idle waiting on others to finish their computations. For each core, if we define these three values as $T_u, T_c$ and $T_s$, respectively, then core utilization can be estimated as $\frac{T_u}{T_u+T_c+T_s}$ \cite{munguia2018alternating}. Figure \ref{fig:utilization} reports average core utilization for each instance. We again see diminishing returns where average utilization decreases as the number of blocks and cores used increases. However, it seems that the average utilization is slightly better as we increase the number of total variables.

\section{Conclusion}
In this paper, we proposed a consensus-based Dantzig-Wolfe decomposition algorithm for loosely coupled large-scale linear programs. As opposed to the standard Dantzig-Wolfe algorithm, we solve the master problem using consensus-based ADMM in a distributed fashion, thus handling circumstances where data is not available centrally and preserving privacy of information of the blocks, and addressed the resulting computational and theoretical challenges. We proved convergence of the algorithm and provided error bounds on the feasibility and optimality gaps. We illustrated our method using an MPI implementation on cutting stock and synthetic instances, and showed that we are able to achieve high accuracy in reasonable time. Although the main objective of our method is to handle decentralized storage of data and privacy concerns, we illustrated potential computational benefits for instances with a large number of variables and in many of the cutting stock instances. To further improve computation time, it is possible to use other algorithms or more sophisticated versions of ADMM to solve the consensus problem. As the difficulty and size of the problems for each block increases, the cost per ADMM iteration can become prohibitive. Certain workarounds involve linearizing the objective of the augmented Lagrangian, yielding computational benefits \cite{boyd}. Other interesting consensus algorithms include a distributed interior point method which might converge faster than first-order distributed methods \cite{bitlisliouglu2017interior}. Finally, as suggested by our experiments, an asynchronous implementation of CDWD has the potential to improve computation times.

\section*{Acknowledgement}
This resesarch has been supported in part by the Office of Naval Research under the Award number N00014-18-1-2075. We would also like to thank George Nemhauser, Natashia Boland and Kibaek Kim for their helpful suggestions on improving the paper.
 
\bibliographystyle{abbrv}
\bibliography{references}
\end{document}